\documentclass[12pt]{article}  
\usepackage{amsthm,amsmath,amsfonts,epsfig,amssymb,latexsym}
\usepackage{mathrsfs}
\usepackage[T1]{fontenc}
\usepackage{graphicx}
\usepackage{enumerate}
\usepackage{xy,xypic}
\usepackage{graphics}
\usepackage{footnote}
\usepackage{scalefnt}
\usepackage{tikz}
\usetikzlibrary{shapes,arrows,fit,calc,positioning}
\usetikzlibrary{matrix}
\xyoption{all}
\title{{\Large The Homotopy Obstructions in Complete Intersections}} 
\author{
 Satya Mandal\footnote{Partially supported by a General Research Grant (no 2301857) from U. of Kansas}
 ~~~and ~~Bibekananda Mishra
 \\ 
{\small University of Kansas, Lawrence, Kansas 66045, USA}\\
{\small {\it  mandal@ku.edu, bibekanandamishra@ku.edu} 
  }\\
 } 
\begin{document}
\renewcommand{\baselinestretch}{1.255}
\setlength{\parskip}{1ex plus0.5ex}
\date{25 December 2016 \\ Revised: 8 August 2017}
\newtheorem{theorem}{Theorem}[section]
\newtheorem{proposition}[theorem]{Proposition}
\newtheorem{lemma}[theorem]{Lemma}
\newtheorem{corollary}[theorem]{Corollary}
\newtheorem{construction}[theorem]{Construction}
\newtheorem{notations}[theorem]{Notations}
\newtheorem{question}[theorem]{Question}
\newtheorem{example}[theorem]{Example}
\newtheorem{definition}[theorem]{Definition} 
\newtheorem{conjecture}[theorem]{Conjecture} 
\newtheorem{remark}[theorem]{Remark} 
\newtheorem{statement}[theorem]{Statement}

\newcommand{\iso}{\stackrel{\sim}{\longrightarrow}}
\newcommand{\sur}{\twoheadrightarrow}
\newcommand{\bD}{\begin{definition}}
\newcommand{\eD}{\end{definition}}
\newcommand{\bP}{\begin{proposition}}
\newcommand{\eP}{\end{proposition}}
\newcommand{\bL}{\begin{lemma}}
\newcommand{\eL}{\end{lemma}}
\newcommand{\bT}{\begin{theorem}}
\newcommand{\eT}{\end{theorem}}
\newcommand{\bC}{\begin{corollary}}
\newcommand{\eC}{\end{corollary}} 
\newcommand{\eop}{\hfill \rule{2mm}{2mm}}
\newcommand{\pf}{\noindent{\bf Proof.~}}
\newcommand{\PD}{\text{proj} \dim}
\newcommand{\lra}{\longrightarrow}
\newcommand{\hra}{\hookrightarrow}
\newcommand{\llra}{\longleftrightarrow}
\newcommand{\Lra}{\Longrightarrow}
\newcommand{\Llra}{\Longleftrightarrow}
\newcommand{\bE}{\begin{enumerate}}
\newcommand{\eE}{\end{enumerate}}
\newcommand{\Sets}{\underline{{\mathrm Sets}}}
\newcommand{\Sch}{\underline{{\mathrm Sch}}}
\newcommand{\ForMe}{\noindent\TCP{{\bf Remarks To Be Removed:~}}}
\newcommand{\pic}{The proof is complete.}
\newcommand{\tcp}{This completes the proof.}

\def\spec#1{\mathrm{Spec}\left(#1\right)}
\def\m{\mathfrak {m}}
\def\CA{\mathcal {A}}
\def\CB{\mathcal {B}}
\def\CP{\mathcal {P}}
\def\CC{\mathcal {C}}
\def\CD{\mathcal {D}}
\def\CE{\mathcal {E}}
\def\CF{\mathcal {F}}
\def\CE{\mathcal {E}}
\def\CG{\mathcal {G}}
\def\CH{\mathcal {H}}
\def\CI{\mathcal {I}}
\def\CJ{\mathcal {J}}
\def\CK{\mathcal {K}}
\def\CL{\mathcal {L}}
\def\CM{\mathcal {M}}
\def\CN{\mathcal {N}}
\def\CO{\mathcal {O}}
\def\CP{\mathcal {P}}
\def\CQ{\mathcal {Q}}
\def\CR{\mathcal {R}}
\def\CS{\mathcal {S}}
\def\CT{\mathcal {T}}
\def\CU{\mathcal {U}}
\def\CV{\mathcal {V}}
\def\CW{\mathcal {W}}
\def\CX{\mathcal {X}}
\def\CY{\mathcal {Y}}
\def\CZ{\mathcal {Z}}

\newcommand{\smallcirc}[1]{\scalebox{#1}{$\circ$}}
\def\BA{\mathbb {A}}
\def\BB{\mathbb {B}}
\def\BC{\mathbb {C}}
\def\BD{\mathbb {D}}
\def\BE{\mathbb {E}}
\def\BF{\mathbb {F}}
\def\BG{\mathbb {G}}
\def\BH{\mathbb {H}}
\def\BI{\mathbb {I}}
\def\BJ{\mathbb {J}}
\def\BK{\mathbb {K}}
\def\BL{\mathbb {L}}
\def\BM{\mathbb {M}}
\def\BN{\mathbb {N}}
\def\BO{\mathbb {O}}
\def\BP{\mathbb {P}}
\def\BQ{\mathbb {Q}}
\def\BR{\mathbb {R}}
\def\BS{\mathbb {S}}
\def\BT{\mathbb {T}}
\def\BU{\mathbb {U}}
\def\BV{\mathbb {V}}
\def\BW{\mathbb {W}}
\def\BX{\mathbb {X}}
\def\BY{\mathbb {Y}}
\def\BZ{\mathbb {Z}}

\newcommand{\TCP}{\textcolor{purple}}
\newcommand{\TCM}{\textcolor{magenta}}
\newcommand{\TCR}{\textcolor{red}}
\newcommand{\TCB}{\textcolor{blue}}
\newcommand{\TCG}{\textcolor{green}}

\def\SA{\mathscr {A}}
\def\SB{\mathscr {B}}
\def\SC{\mathscr {C}}
\def\SD{\mathscr {D}}
\def\SE{\mathscr {E}}
\def\SF{\mathscr {F}}
\def\SG{\mathscr {G}}
\def\SH{\mathscr {H}}
\def\SI{\mathscr {I}}
\def\SJ{\mathscr {J}}
\def\SK{\mathscr {K}}
\def\SL{\mathscr {L}}
\def\SN{\mathscr {N}}
\def\SO{\mathscr {O}}
\def\SP{\mathscr {P}}
\def\SQ{\mathscr {Q}}
\def\SR{\mathscr {R}}
\def\SS{\mathscr {S}}
\def\ST{\mathscr {T}}
\def\SU{\mathscr {U}}
\def\SV{\mathscr {V}}
\def\SW{\mathscr {W}}
\def\SX{\mathscr {X}}
\def\SY{\mathscr {Y}}
\def\SZ{\mathscr {Z}}

\def\bfA{{\bf A}}
\def\bfB{{\bf B}} 
\def\bfC{{\bf C}} 
\def\bfD{{\bf D}} 
\def\bfE{{\bf E}} 
\def\bfF{{\bf F}} 
\def\bfG{{\bf G}} 
\def\bfH{{\bf H}} 
\def\bfI{{\bf I}} 
\def\bfJ{{\bf J}} 
\def\bfK{{\bf K}} 
\def\bfL{{\bf L}} 
\def\bfM{{\bf M}} 
\def\bfN{{\bf N}} 
\def\bfO{{\bf O}} 
\def\bfP{{\bf P}} 
\def\bfQ{{\bf Q}} 
\def\bfR{{\bf R}} 
\def\bfS{{\bf S}} 
\def\bfT{{\bf T}} 
\def\bfU{{\bf U}} 
\def\bfV{{\bf V}} 
\def\bfW{{\bf W}} 
\def\bfX{{\bf X}} 
\def\bfY{{\bf Y}} 
\def\bfZ{{\bf Z}} 

\maketitle

\section{Introduction}

For (smooth) affine schemes $X=\spec{A}$ over algebraically closed fields $k$ and for projective $A$-modules $P$,  with $\dim X=d=rank(P)$, 
N. Mohan Kumar and M. P. Murthy considered top Chern classes $C^d(P)\in CH^d(X)$, in the Chow Group of zero cycles,
as obstruction for $P$ to split off a free direct summand (\cite{MoM, Mk1, Mk2,  Mu, MMu}). In deed, the results in \cite{Mu, MMu}, were
fairly finalistic.  
For an ideal $I$, an obstruction class $\zeta(P, I)\in CH^d(X)$ was written down \cite{MMu}.
%
It was established that  there is a surjective map $P\sur I$ if and only if $\zeta(P, I)=0$ and there is a surjective map $P\sur \frac{I}{I^2}$.

Subsequent to that,  based on some Homotopy Relations (see Lemma \ref{noriPizero}), Madhav
V. Nori (around 1990) laid out a set of ideas to deal with the questions of such obstructions in  broader contexts, like when $X$ is a regular or  a noetherian affine scheme.
 These were communicated verbally to some 
in a very informal and open ended manner. Because of the nature of these communications, not everyone heard the same thing and
 these ideas took the form of some folklores. 
 As a result, 
 versions of this set of ideas available  (or not) in the literature (e. g. \cite{M3, MV, MS,BS1,BS2, BS3, BK})
 have been up to the interpretations and 
adaptations by  the  recipients of these communications, much to their credit, and the stated hypotheses may differ.
 Because of the openendedness and broadness of these ideas, they appeared to be more like a research Program ({\it the Homotopy Program})
   to some, 
which is how   we would refer  to the same in this article.
Sometimes it may even be difficult to say whether certain part of the program was actually explicitly articulated by Nori 
or  were part of the adaptations by others.
There is no  systematic exposition of this program available in the literature and 
 certain aspects  failed to receive deserved
traction. Nori never classified these  as conjectures or otherwise.
However, some results followed too quickly \cite{M3},
to treat them as anything less than conjectures, subject to further fine tuning.
%

Analogy to the Obstruction Theory for vector bundles (\cite{St})  was the main backdrop behind this program
and  central to this Program was the Homotopy conjecture of Nori. The following is the statement of the Homotopy Conjecture from \cite{M3},
which would most likely be an adaptation by the respective author.

\begin{conjecture}[Homotopy Conjecture]\label{homoConj}
Suppose $X=\spec{A}$ is a smooth affine variety, with $\dim X=d$. Let $P$ be a projective $A$-module of rank $r$ and $f_0:P\sur I$ 
be  a surjective 
homomorphism, onto an ideal $I$ of $A$. Assume $Y=V(I)$ is smooth with $\dim Y=d-r$.
Also suppose $Z=V(J)\subseteq \spec{A[T]}=X\times \BA^1$ is a smooth subscheme, such that $Z$ intersects $X\times 0$ transversally in $Y\times 0$.
Now, suppose that $\varphi: P[T] \sur \frac{J}{J^2}$ is a surjective map such that $\varphi_{|T=0}=f_0\otimes \frac{A}{I}$.  Then, there is a
surjective 
map $F:P[T]\sur J$ such that (i) $F_{|T=0}=f_0$ and (ii) $F_{|Z}=\varphi$.
\end{conjecture}

The best result, up to date, on this Conjecture \ref{homoConj} is due to Bhatwadekar and Keshari \cite{BK}.
While the Conjecture  \ref{homoConj} would fail without the regularity  hypothesis   \cite[Example 6.4]{BS1},
existing results (see \cite{M3, BS1, BK}) indicate that with suitable  hypotheses the regularity  and/or transversality hypotheses may be spared.
However, the Conjecture  \ref{homoConj}, as stated,  would  even fail in some cases when $A$ is regular (see \cite[Example 3.15]{BS1}), which are, conjecturally, the exceptions.
In analogy to the obstructions $\zeta(P, I)\in CH^{d}(X)$ mentioned above (\cite{MMu}), the 
objective of the Homotopy Conjecture \ref{homoConj} was to detect, for an ideal $I$, when  a surjective map 
$f:P \sur \frac{I}{I^2}$  would lift to a surjective 
map $F: P\sur I$, using the homotopy relations given by surjective maps $\varphi: P[T] \sur \frac{J}{J^2}$, where $J$ is an ideal in $A[T]$ 
(see Lemma \ref{noriPizero}).
%

In this article, we will discuss the Homotopy Program, only in the complete intersections case, that is when $P=A^n$ is free, with $n\geq 2$. 
To clarify the Homotopy Obstruction set of Nori, in this complete intersection case, 
let ${\mathcal LO}(A, n)$ denote the set of all pairs $(I, \omega)$, where $I$ is an ideal of $A$ and $\omega: A^n \sur \frac{I}{I^2}$
is a surjective homomorphism. By substituting $T=0, 1$, we obtain two maps $\diagram {\mathcal LO}(A, n) & {\mathcal LO}(A[T], n)\ar[l]_{T=0}\ar[r]^{T=1} &
{\mathcal LO}(A, n)\enddiagram$.
This leads to an equivalence relation (using chains) on ${\mathcal LO}(A, n)$ and a set of all equivalence classes $\pi_0({\mathcal LO}(A,n))$.
In the recent past, a similar Homotopy Obstruction set (pre-sheaf),
$\pi_0(Q_{2n})(A)$ was considered in \cite{F}, to serve the same purpose. We clarify (see Lemma \ref{noriPizero}) that $\pi_0(Q_{2n})(A)$  coincides with $\pi_0({\mathcal LO}(A,n))$.
This puts some of the developments \cite{F, M2} ({\it now retracted})
in the recent past, in the framework of the Homotopy Program of Nori.
In  deed, in this article, we mostly investigate the structure of this obstruction set $\pi_0({\mathcal LO}(A,n))=\pi_0(Q_{2n})(A)$.

First, out of necessity, we prove a quadratic version of Lindel's Theorem (\cite{L}, Bass-Quillen Conjecture), 
on extendibility of projective modules $P$ over polynomial rings $R=A[T]$, where $A$ is a regular ring containing a field, as follows.
\bT\label{IntroorthPopescu} 
Suppose $A$ is a regular ring over a field $k$, with
$1/2\in k$ and $\dim A=d$. Let $R=A[T_1, \ldots, T_n]$ be a polynomial ring. 
Suppose  $q$ is an isotropic quadratic form over $k$, with $rank(q)=r$. Suppose 
$(P, \varphi)$ is a quadratic space over $R$. Write $(\overline{P}, \overline{\varphi}):=(P, \varphi)\otimes \frac{R}{(T_1, \ldots, T_n)R}$.
Assume $(\overline{P}, \overline{\varphi})$ is locally trivial {\rm (in the sense clarified in (\ref{orthLindel})).}
Then, $(P, \varphi)$ is extended from $A$.  
\eT
While Theorem  \ref{IntroorthPopescu}  is significant by its own right, the theorem completely removes the "infinite field" condition,  from the methods in   \cite{F, M2}.

For commutative noetherian rings $A$,  we establish (see \S  \ref{secInvolution}) the existence of a natural involution on $\pi_0({\mathcal LO}(A,n))$, for all $n\geq 2$,
which is a key tool in this article. 
The possibility of a group structure on the obstruction set $\pi_0({\mathcal LO}(A,n))$ has always been a part of the Homotopy Program, particularly,
 in the upper half of the range of $n$. 
 We establish a group structure on $\pi_0({\mathcal LO}(A,n))$ as follows.  

\bT\label{INTROabelianGroup}
Suppose $A$ is a regular ring over a field $k$, with $1/2\in k$, with $\dim A=d$. 
Let $n\geq 2$ be an integer, with $2n\geq d+2$.
Then, there is a structure of an abelian group on $\pi_0({\mathcal LO}(A, n))=\pi_0(Q_{2n})(A)$. The addition is 
determined as follows: Suppose $x=\zeta(K, \omega_K), y=\zeta(I, \omega_I)\in \pi_0(Q_{2n})(A)$,
where $\zeta: {\mathcal LO}(A, n)\lra \pi_0(Q_{2n})(A)$ is the natural map and 
 $(K, \omega_K), (I, \omega_I)\in {\mathcal LO}(A, n)$, such that
 $height(K)\geq n$ and $K+I=A$. Then,
$$
x+y=\zeta(KI, \omega_K\star\omega_I) \in \pi_0(Q_{2n})(A)
$$
where $ \omega_K\star\omega_I: A^n \sur \frac{KI}{(KI)^2}$ is the unique surjective map determined by 
$\omega_K$ and $\omega_I$.
\eT

In deed, another definition of  an Obstruction group $E^{d}(A)$, where $d=\dim A$, of zero cycles 
 was outlined by Nori, by considering the free abelian group generated by $\{(m, \omega)\in {\mathcal LO}(A, \dim A): m\in \max(A)\}$, and the relations 
 were obtained 
 using homotopy equivalences ({\it see, for example, \cite[\S 4, pp.175-176]{BS1}, where it was denoted by $E(A)$}). 
 By suitable adaptations, 
 Bhatwadekar and Sridharan \cite{BS2} defined  obstruction groups $E^n(A)$, for each co-dimension $n\geq 0$, 
 where $A$ was assumed to be any noetherian commutative ring (see \S~\ref{eulerDefSec}).
 These groups $E^n(A)$ are known as Euler class groups.
  Ever since, relationship between the obstructions $\pi_0({\mathcal LO}(A, n))$ 
  and $E^n(A)$
 remained an open question in the Homotopy Program, which is settled follows. 

\bT\label{IntrogrMappi02En}
 Suppose $A$ is a regular ring over a field $k$, with $1/2\in k$, with $\dim A=d$. 
Assume $n$ is an integer such that $n\geq 2$ and $2n\geq d+2$. Then, there is a surjective homomorphism, 
$\rho: E^n(A) \sur \pi_0({\mathcal LO}(A, n))$. Further, this homomorphism is an isomorphism, if for orientations $(I, \omega_I)\in {\mathcal LO}(A, n)$,
with $height(I)\geq n$,
its triviality in $\pi_0({\mathcal LO}(A, n))$ implies $\omega_I$ lifts to a surjective map $A^n\sur I$.

In particular, by \cite{BK},   $\rho$ is an isomorphisms, if $A$ is essentially smooth and $k$ is an  infinite perfect field,
and $2n\geq d+3$ (see Theorem \ref{useBKNov14}).
 \eT
\begin{remark}\label{pfIncomplete}{\rm
After  this article was posted in arXiv, we became aware 
of the   article \cite{AF}, which has some overlap with the results in this article. The results in this article would be an improvement upon those in \cite{AF}, and the
 methods  are completely different. Our main results hold for any regular ring $A$ containing a field $k$, with $1/2\in k$, while main results in \cite{AF} are valid for smooth rings  over  infinite perfect fields $k$, with $1/2\in k$. Due to uncertainty  regarding the validity of \cite[Theorem 3.2.8]{F}, ({\it now retracted} \cite{F2}), 
there is no valid proof of the injectivity of the map $s$ in Theorem 3.1.9 in \cite{AF}, when $2n=d+2$. In this article (Theorem \ref{useBKNov14}), we use \cite[Theorem 4.13]{BK} to establish the injectivity, when $2n\geq d+3$. 
}
\end{remark}
\begin{remark}{\rm 
Let $\widetilde{CH}^d(A)$ denote the Chow-Witt group, defined by Barge and Morel \cite{BM}. 
It is known that, under the hypotheses of (\ref{IntrogrMappi02En}), \\
$\pi_0({\mathcal LO}(A, d)) \cong \pi_0\left(Q_{2d}(A) \right) \cong \widetilde{CH}^d(A)$,
where the first isomorphism is established in Lemma \ref{noriPizero} and the second isomorphism was studied in \cite{AF2}. 

It follows from Theorem \ref{IntrogrMappi02En} that, under the same hypotheses, 
there is an isomorphism $E^d(A)\iso \widetilde{CH}^d(A)$. This settles a long standing open question (see \cite[Remark 1.33(2)]{Mor}),
as was also claimed  in \cite[Theorem 3 (3)]{AF}. However, the proof of \cite[Theorem 3 (3)]{AF} remains incomplete, as clarified above (\ref{pfIncomplete}). 
}
\end{remark}


The organization of this article is fairly sequential, as the results are described above. 
We add that, in Section \ref{secHomotopy}, 
 some important homotopy theorems were established. This, in particular,  establishes that 
 homotopy relations described above, is actually an equivalence relation, when the ring $A$ is regular over a field $k$, with $1/2\in$.
In Section \ref{eulerDefSec}, we give some preliminaries about Euler class groups. In Section \ref{EOsection}, we provide some background 
regarding Elementary Orthogonal Subgroups ${\SE}O(A, q_{2n+1})$.
 
\noindent{\bf Acknowledgement:} {\it 
Thanks are  due to Marco Schlichting for his continued academic 
support and for pointing to the paper of M. Ojanguren \cite{O1}. We also thank the referee for careful reviewing and for pointing to the reference to \cite{Rr}.}

\section{Notations and Preliminaries} \label{prelimSec}
First, we establish some notations, some of which may be standard.
\begin{notations}\label{nota}{\rm 
Throughout, $k$ will denote a field (or a ring),
with $1/2\in k$ and $A$ will denote a commutative noetherian rings.
Denote\\
 $ \tilde{q}_{2n+1}=\sum_{i=1}^nX_iY_i+Z(Z-1)$,
\begin{equation}
 \label{Q2nSept18}
\SA_{2n}=\frac{k[X_1, \ldots, X_n, Y_1, \ldots, Y_n, Z]}{\left(\tilde{q}_{2n+1}\right)} \quad {\rm and}\quad
 Q_{2n}=\spec{\SA_{2n}}.
 \end{equation}
Accordingly, for a commutative ring $A$, denote 
$$
 Q_{2n}(A)=\left\{(s; f_1,\ldots, f_n; g_1, \ldots, g_n)\in A^{2n+1}: \sum_{i=1}^nf_ig_i+s(s-1)=0 \right\}
 $$
For ${\bf v}=(s; f_1, \ldots, f_n; g_1, \ldots, g_n)\in Q_{2n}(A)$, denote the ideals 
 $$
 \left\{
 \begin{array}{l}
 \BI({\bf v}):=
 (f_1, \ldots, f_n, s)A\\
 \BJ({\bf v}):=
 (f_1, \ldots, f_n, 1-s)A\\
 \end{array}
 \right.
~~{\rm and~surjective~maps}~~
 \left\{
 \begin{array}{l}
\omega_{\bf v}:A^n\to \frac{\BI({\bf v})}{\BI({\bf v})^2}\\
\omega_{\bf v}':A^n\to \frac{\BJ({\bf v})}{\BJ({\bf v})^2}\\
 \end{array}
 \right.
 $$
%
 %
  Also, $\Sets$ will denote the category of sets.
 The homotopy pre-sheaves are given by the pushout diagrams in $ \Sets$:
 \begin{equation}\label{DiaPushpoi0Q2n}
 \diagram 
  Q_{2n}(A[T]) \ar[r]^{T=0}\ar[d]_{T=1} & Q_{2n}(A)\ar[d]\\
Q_{2n}(A) \ar[r] & \pi_0\left(Q_{2n}\right)(A)\\ 
\enddiagram 
\end{equation}
Further, consider the quadratic form
 $q_{2n+1}=\sum_{i=1}^nX_iY_i+Z^2$, and denote
\begin{equation}
 \label{Q2nSept236teen}
 \SB_{2n}=\frac{k[X_1, \ldots, X_n, Y_1, \ldots, Y_n, Z]}{\left(q_{2n+1}-1\right)},\qquad  Q_{2n}'=\spec{\SB_{2n}}.
 \end{equation}
 Accordingly,
 $$
 Q_{2n}'(A)=\left\{(s; f_1,\ldots, f_n;  g_1, \ldots, g_n)\in A^{2n+1}: \sum_{i=1}^nf_ig_i+s^2-1=0 \right\}
 $$
 Since, $1/2\in k$ there is an isomorphism  $\alpha: \SA_{2n}\iso \SB_{2n}$, which induces bijective correspondence
 $$
 \alpha: Q_{2n}'(A) \iso Q_{2n}(A) \quad{\rm and~its~inverse}\quad  \beta: Q_{2n}(A) \iso Q_{2n}'(A)
 $$ 
 Now, the action of the  Orthogonal groups $O(A, q_{2n+1})$  on $Q_{2n}'(A)$ translates to an action of 
 $O(A, q_{2n+1})$  on $Q_{2n}(A)$ as follows:
 $$
\forall~{\bf v}\in Q_{2n}(A), M\in O\left(A,q_{2n+1}\right) ~~{\rm define} ~~{\bf v}*M:=\alpha\left(\beta({\bf v})M \right)
$$

}
\end{notations} 
The  local orientations of an ideal 
are defined as follows.
\bD{\rm
Suppose $A$ is a commutative  ring and $I$ is an ideal in $A$. For an integer $n\geq 2$, 
a \TCP{local $n$-orientation} of $I$ is a pair $(I, \omega)$, where
  $\omega:A^n \sur I/I^2$ is a a surjective homomorphism.
Such a local $n$-orientation is determined by any set of elements $f_1, \ldots, f_n\in I$ such that $I =(f_1, \ldots, f_n)+I^2$.  Given such a set of generators $f_1, \ldots, f_n$ of $I/I^2$, there is an element $s\in I$ and $g_1, \ldots, g_n$ such that  $\sum_{i=1}^nf_ig_i+s(s-1)=0$. Note,
$$
(s; f_1, \ldots, f_n; g_1, \ldots, g_n)\in Q_{2n}(A).
$$
$$
{\rm Write}\quad 
\zeta(I, \omega):=[(s; f_1, \ldots, f_n; g_1, \ldots, g_n)]\in \pi_0\left(Q_{2n}(A)\right)
$$
This association is well defined (\cite[Theorem 2.0.7]{F}). We refer to $\zeta(I, \omega)$,
as the 
\TCP{homotopy obstruction class}. 
The set of all $n$-orientations $(I, \omega_I)$ will be denoted by ${\mathcal LO}(A, n)$. %
 Therefore, we have a commutative diagram
\begin{equation}\label{zetaDiag}
\diagram
Q_{2n}(A)\ar[rd]^{\zeta_0}\ar[d]_{\eta} & \\
{\mathcal LO}(A, n) \ar[r]_{\zeta} &  \pi_0\left(Q_{2n}(A)\right)\\
\enddiagram
\qquad{\rm where}\quad \eta({\bf v})=\left(\BI({\bf v}), \omega_{{\bf v}}\right).
\end{equation}
%
Analogous to the definition of $\pi_0(Q_{2n})(A)$, we define $\pi_0(\CL\CO(A, n))$, by the pushout diagram
\begin{equation}\label{DiaPushpiZeroOAn}
 \diagram 
{\mathcal LO}(A[T], n) \ar[r]^{T=0}\ar[d]_{T=1} &  {\mathcal LO}(A, n)\ar[d]\\
{\mathcal LO}(A, n) \ar[r] & \pi_0({\mathcal LO}(A,n))\\ 
\enddiagram 
\qquad {\rm in}\quad \underline{Sets}.
\end{equation}
}
\eD
%
The above homotopy  obstruction set  $\pi_0({\mathcal LO}(A,n))$ was among the ideas envisioned by Nori (around 1990).
Following lemma establishes that $\pi_0(Q_{2n})(A)$ is in bijection with $\pi_0({\mathcal LO}(A,n))$.

\bL\label{noriPizero}
Suppose $A$ is a commutative noetherian ring, with $\dim A=d$ and $n\geq 2$ is an integer.
Then, the map $\zeta:  {\mathcal LO}(A, n) \lra \pi_0(Q_{2n})(A)$ induces bijective map
$$
\overline{\zeta}: \pi_0({\mathcal LO}(A,n))\iso \pi_0(Q_{2n})(A).
$$
The inverse  map is induced by $\eta$.
\eL
\pf Given a $\tilde{H}\in  \pi_0({\mathcal LO}(A[T],n))$, $\eta(H(T))= \tilde{H}$ for some $H(T)\in Q_{2n}(A[T])$. 
So, $\zeta(\tilde{H}(0))= \zeta_0(H(0))= \zeta_0(H(1))= \zeta(\tilde{H}(1))$.
This establishes that $\zeta$ factors through a set theoretic map $\overline{\zeta}: \pi_0({\mathcal LO}(A,n))\lra \pi_0(Q_{2n})(A)$.
Since $\eta$ is on to , so is $\overline{\zeta}$. There is also a well defined map 
$\overline{\eta}: \pi_0(Q_{2n})(A) \lra \pi_0({\mathcal LO}(A,n))$ induced by $\eta$. It is clear that $\overline{\eta}\overline{\zeta}=1$.
So, $\overline{\zeta}$ is also one to one and hence
 is bijective.
\pic $\eop$

The following  "moving lemma argument" is fairly standard. A number of  variations of the same (\ref{movingLemm}) 
would be among the frequently used tools for the rest of our discussions.
\bL[Moving Lemma]\label{movingLemm}
Suppose $A$ is a commutative noetherian ring with $\dim A=d$ and $n\geq 2$ is an integer such that $2n\geq d+1$.
Let $K\subseteq A$ be an ideal with $height(K)\geq n$ and $(I, \omega_I)\in \CL\CO(A, n)$. 
Then, there is an element ${\bf v}=(s; f_1, \ldots, f_n; g_1, \ldots, g_n)\in Q_{2n}(A)$ such that 
$\eta({\bf v}) =(I, \omega_I)$. Further, with $J=\BJ({\bf v})$, we have $height(J)\geq n$ and $J+K=A$.
\eL
\pf We use the standard basis $e_1, \ldots, e_n$ of $A^n$. 
Let $a_1, \ldots, a_n\in I$ be such that $\omega_I(e_i)=a_i+I^2$. 
So, $I=(a_1, \ldots, a_n)+I^2$. Using Nakayama's Lemma, 
there is an element $t\in I$,
such that $t(1-t)=\sum_{i=0}^na_ib_i$ for  some $b_1, \ldots, b_n\in A$ and $I=(a_1, \ldots, a_n, t)$. 
({\it Readers are referred to \cite{M1} regarding generalities on Basic Element Theory and generalized
dimension functions}.)
Write 
$$
\SP=\left\{\wp\in \spec{A}: t\notin \wp, ~{\rm and~either}~K\subseteq \wp~{\rm or}~height(\wp)\leq n-1 \right\}
$$
There is a generalized dimension function (see \cite{M1})  $\delta: \SP \lra \BN$, such that $\delta(\wp)\leq n-1~\forall~\wp\in \SP$.
Now $(a_1, \ldots, a_n, t^2)\in A^{n+1}$ is basic on $\SP$. So, there are $\lambda_1, \ldots, \lambda_n\in A$ such that 
$(a_1+\lambda_1t^2, \ldots, a_n+\lambda_nt^2)\in A^{n}$  is basic on $\SP$. For $i=1, \ldots, n$, denote $f_i:=a_n+\lambda_it^2$.
Then, $\omega_I(e_i)=f_i+I^2$ and hence $I=(f_1, \ldots, f_n)+I^2$. 
Using Nakayama's Lemma, there is an element $s\in I$ such that $(1-s)J\subseteq (f_1, \ldots, f_n)$.
Hence, $s(1-s)=\sum_{i=1}^nf_ig_i$ for some $g_1, \ldots, g_n$. Now, 
${\bf v}=(s; f_1, \ldots, f_n; g_1, \ldots, g_n)\in Q_{2n}(A)$ has the desired properties. 
\pic $\eop$

\begin{remark}{\rm 
 In Moving Lemma \ref{movingLemm}, we proved that, given \\
 ${\bf u}=(t; a_1,\ldots, a_n; b_1, \ldots, b_n) \in Q_{2n}(A)$,
there is\\
 ${\bf v}=(s; f_1, \ldots, f_n; g_1, \ldots, g_n)\in Q_{2n}(A)$, such that 
\bE
\item  $\forall~i=1, \ldots, n~f_i=a_i+\lambda_it^2$,
\item $\eta({\bf u})=\eta({\bf v})$,
\item $height(\BJ({\bf v}))\geq n$ and $\BJ({\bf v})+K=A$.
\eE
}\end{remark}
For the convenience of our discussions, we include some auxiliary notations.
\begin{notations}\label{auxNota}{\rm 
Throughout, $A$ will denote a noetherian commutative ring, with $\dim A=d$.
\bE
\item For an $A$-module $M$ and an ideal $I$ of $A$,  homomorphisms  $f: M\lra \frac{I}{I^2}$ would be identified with the induced 
maps $\frac{M}{IM} \lra \frac{I}{I^2}$.
\item Let $I_1, I_2$ be two ideals, with $I_1+I_2=A$. For an integer $n\geq 2$, for $i=1, 2$ let $\omega_i:A^n\sur \frac{I_i}{I_i^2}$
be two surjective maps. Then, $\omega_1\star\omega_2: A^n\sur \frac{I_1I_2}{(I_1I_2)^2}$ will denote the unique surjective map 
induced by $\omega_1, \omega_2$. 
\item For integers $n\geq 2$, denote ${\bf 0}:=(0; 0, \ldots, 0; 0, \ldots, 0)\in Q_{2n}(A)$
and  ${\bf 1}:=(1; 0, \ldots 0; 0, \ldots, 0)\in Q_{2n}(A)$. Either one of them could be a candidate for the base point of $Q_{2n}(A)$.
\item An element $H(T)\in  Q_{2n}(A[T])$ or $H(T)\in  Q_{2n}'(A[T])$ would be referred to as a homotopy.
\item We caution the readers that notations in this section would be part of our standard notations, throughout this article. In particular, that 
would include $\BI({\bf v})$, $\BJ({\bf v})$, $\omega_{{\bf v}}$, $\eta$, $\zeta_0$, $\zeta$ and others.
\eE
}
\end{notations}



\section{A Quadratic Version of Lindel's Theorem}\label{Lindel}

For an essentially smooth ring $A$ over a field $k$, and a polynomial ring $R=A[T_1, \ldots, T_n]$,
Lindel \cite{L} settled Bass-Quillen Conjecture, by proving that  finitely generated projective  $R$-modules are extended from $A$. 
By the Desingularization Theorem of Popescu (\cite{P}, \cite[Corollary 1.2]{Sw}), it follows that the same is true, when $A$ is  any regular ring 
containing a field $k$ (see \cite[Theorem 2.1]{Sw}). In this section, we give a version of the same for Quadratic spaces.

For the convenience of the readers, we recall the
following definition.
\bD
Suppose $A$ is a noetherian commutative ring, with $1/2\in A$. Then, a Quadratic space on $A$ is a pair $(P, \varphi)$, where $P$ is a 
finitely generated projective $A$-module and $\varphi: P \iso P^*$ is an isomorphism, such that $\varphi^*=\varphi$.
Maps between two quadratic spaces over $A$, are called orthogonal maps and such isomorphisms are called isometries.
\eD

First, we prove the Quadratic analogue of Lindel's theorem for smooth rings $A$ over perfect fields, as follows.

\bP\label{orthLindel} 
Suppose $A$ is an essentially smooth  ring over a  prefect field $k$, with
$1/2\in k$ and $\dim A=d$. Let $R=A[T_1, \ldots, T_n]$ be a polynomial ring. 
Suppose  $q$ is an isotropic quadratic form over $k$, with $rank(q)=r$. Suppose 
$(P, \varphi)$ is a quadratic space over $R$. Write $(\overline{P}, \overline{\varphi}):=(P, \varphi)\otimes \frac{R}{(T_1, \ldots, T_n)R}$.
Assume $(\overline{P}, \overline{\varphi})$ is locally isometric to $(A^n, q)$; meaning $\forall~\wp\in \spec{A}$,  $(\overline{P}, \overline{\varphi})_{\wp} \cong 
(A^n, q)_{\wp}$ are isometric. {\rm (This hypothesis is referred to as "local triviality property").}
Then, $(P, \varphi)$ is extended from $A$.  
\eP
\pf 
We prove by induction on $\dim A=d$.
Suppose $d=\dim A=0$. Then, $A=k$ is a field. By hypothesis $WittIndex(\overline{P})\geq 1$.
So,  the theorem is valid, by the Theorem of Ojanguren \cite{O1} (also see \cite[pp.425, Thm 6.2.6]{K}). 
So, we assume $\dim A\geq 1$. 
We can also assume that $A$ is local (see \cite[pp. 419, Thm 5.3.4]{K}). 
By hypothesis on local triviality of $(\overline{P}, \overline{\varphi})$, there is an isometry 
\TCP{$\sigma_0:(A^r, q) \iso (\overline{P}, \overline{\varphi})$}.
%
By  Theorem of Lindel \cite{L} (see \cite[7.1.1]{M1}), there is local subring $B\subseteq A$ such that
\bE
\item $B=k[X_1, \ldots, X_d]_M$ where $X_1, \ldots, X_n$ are variables, $f(X_1)\in k[X_1]$ and
 $M=(f(X_1), X_2, \ldots, X_d)$ is a maximal ideal.
\item There is an element $h\in MB$ such that the inclusion map
$B\hra A$ is an \TCP{analytic isomorphism}; meaning 
\begin{equation}\label{anaISO}
\left\{
\begin{array}{l}
A=B+Ah\\
\forall~n\in \BN\quad Bh^n=B\cap Ah^n.\\
\end{array}
\right.
\quad{\rm Consequently,}\quad 
\diagram
B \ar[r]\ar[d] & A\ar[d]\\
B_h \ar[r] & A_h\\
\enddiagram 
\end{equation}
is a patching  diagram (see \cite{O2, R} regarding such patching diagrams). 
\eE
We  denote 
$$
F'=B^n,~~ F'[\underline{T}]=F'\otimes B[T_1, \ldots, T_n],  ~~F=A^n,~~ F[\underline{T}]=F\otimes A[T_1, \ldots, T_n].
$$
So, \TCP{$\sigma_0: (F, q) \iso (\overline{P}, \overline{\varphi})$}, above, is an isometry.
Since, $\dim A_h\leq d-1$, by induction, there is a quadratic space $(Q, \psi)$ on $A_h$, and an isometry
$$
\sigma_1:(Q, \psi)\otimes A_h[T_1, \ldots, T_n] \iso (P_h, \varphi_h).
$$
Let "overline" denote modulo $(T_1, \ldots, T_n)$.
Then $\sigma$ induces an isometry,
$$
\overline{\sigma}_1: (Q, \psi) \iso (\overline{P}_h, \overline{\varphi}_h). 
\quad{\rm So,}~
\tau:=\bar{\sigma}_1^{-1}(\sigma_0)_h: (F_h, q) \iso (Q, \psi)  ~{\rm is~an~isometry.}
$$
 Further, consider the obvious isometry \TCP{$\sigma_2: (F'\otimes A_h, q) \iso (F_h, q)$}.
Combining, all these, we have commutative diagram of isometries 
$$
\diagram
(F'\otimes A_h[\underline{T}], q) \ar[r]^{~~~\sigma_2\otimes 1}\ar@{-->}[drr]_{~~\sigma} & (F_h[\underline{T}], q)\ar[r]^{\tau\otimes 1} & (Q[\underline{T}], \psi\otimes 1)\ar[d]^{\sigma_1}\\
& &P_h\\
\enddiagram
$$
where $\sigma$ is defined by composition. Notice \TCP{$\sigma\otimes \frac{A_h[\underline{T}]}{(\underline{T})}=(\sigma_0)_h\sigma_2$}.
With respect to the patching diagram 
(\ref{anaISO}), consider the patching  diagram (see \cite{O2, R})
$$
\diagram
(\widetilde{P}, \widetilde{\varphi}) \ar[rr]\ar[d] && (P, \varphi)\ar[d]\\
(F'_h[\underline{T}], q) \ar[r] & (F'\otimes A_h[\underline{T}], q) \ar[r]_{\qquad\sigma} & (P, \varphi)_h\\
\enddiagram 
$$
Here $\widetilde{P}$ is obtained by this patching $F'_h[\underline{T}]$ and $P$,  via $\sigma$. 
Then, $\widetilde{P}$ has a structure of a quadratic space $(\widetilde{P}, \widetilde{\varphi})$ \cite[Theorem 8]{O2},
making the above a patching diagram of quadratic spaces. Since, $\sigma\otimes \frac{A_h[\underline{T}]}{(\underline{T})}=(\sigma_0)_h\sigma_2$,
it follows that, there is an isometry $(F', q) \iso (\widetilde{P}, \widetilde{\varphi})\otimes \frac{B[\underline{T}]}{(\underline{T})}$.
(i.e. local triviality property is preserved). 
%
%
Now, replacing $A$ by $B$ and $P$ by $\widetilde{P}$, we can assume $A=k[X_1, \ldots, X_d]_M$ with $M=(f(X_1), X_2, \ldots, X_d)$.

Write $A_0=k[X_1, \ldots, X_{d-1}]_{\m}$, where $\m=(f(X_1), X_2, \ldots, X_{d-1})$. Then, $A_0[X_d] \hra A$
is an \TCP{analytic isomorphism} along $X_d$. Now, the new patching diagram looks like
$$
\diagram 
A_0[X_d] \ar[r] \ar[d] & A\ar[d]\\
A_0[X_d, X_d^{-1}] \ar[r] & A_{X_d}\\
\enddiagram
$$
Repeating the same method, 
$$
(P, \varphi) \cong \left(\widetilde{P}, \widetilde{\varphi}\right) \otimes_{A_0[X_d, T_1, \ldots, T_n]} A[T_1, \ldots, T_n]
$$
 where  $\left(\widetilde{P}, \widetilde{\varphi}\right)$  is a quadratic space over 
$A_0[X_d, T_1, \ldots, T_n]$. 
Since $\dim A_0=d-1$, by induction, $\left(\widetilde{P}, \widetilde{\varphi}\right)$ is extended from $A_0$,
and hence from $A_0[X_d]$. Therefore, $(P, \varphi)$ is extended from $A$.
\tcp $\eop$


\begin{remark}{\rm 
Our interest in this question of extendibility of Quadratic spaces over polynomial rings was triggered,
by the reference to \cite[Theorem 3.3.7]{AHW} in the proof of \cite[Theorem 1.0.6]{F}. It was confusing that a natural quadratic 
analogue of Lindel's theorem \cite{L} would not be considered in \cite{F}, instead, if such was available in the literature.
The statement of  \cite[Theorem 3.3.6]{AHW} obscures the simplicity of such an analogue. 
 However, the referee points to the references \cite{Pr1, Pr2} and most importanly \cite[Theorem 1.3]{Rr}. 
 The result of Ravi Rao \cite[Proposition 1.3]{Rr} is, perhaps, the most significant result on this, available in the literature, which 
 can be viewed as a local case of Proposition \ref{orthLindel}, above.
 Proofs of both \cite[Proposition 1.3]{Rr} and Proposition \ref{orthLindel} would be similar to 
the proof of Lindel's theorem \cite{L}. We  further strengthen Proposition \ref{orthLindel}, as follows (\ref{orthPopescu}), by using Popescu's Desingularization Theorem \cite{P}.
}
\end{remark}

 For the convenience of our discussions, we state the following lemma, which can be checked locally.

\bL\label{HomPQ}
Let $R$ be a noetherian commutative ring and  $A \subseteq R$ be a noetherian subring. 
Let $P, Q$ be two finitely generated projective $A$-modules. 
Then, the map $Hom_A(P, Q) \lra Hom_R(P\otimes R, Q\otimes R)$ is injective.
\eL 
%


Now, we use Popescu's Desingularization Theorem (\cite{P}, \cite[Corollary 1.2]{Sw}),
to remove the perfectness condition in (\ref{orthLindel}), as follows (\ref{orthPopescu}).

\bT\label{orthPopescu} 
Suppose $A$ is a regular ring over a field $k$, with
$1/2\in k$ and $\dim A=d$. Let $R=A[T_1, \ldots, T_n]$ be a polynomial ring. 
Suppose  $q$ is an isotropic quadratic form over $k$, with $rank(q)=r$. Suppose 
$(P, \varphi)$ is a quadratic space over $R$. Write $(\overline{P}, \overline{\varphi}):=(P, \varphi)\otimes \frac{R}{(T_1, \ldots, T_n)R}$.
Assume $(\overline{P}, \overline{\varphi})$ is locally trivial {\rm (in the sense clarified in (\ref{orthLindel})).}
Then, $(P, \varphi)$ is extended from $A$.  
\eT
\pf   For any matrix $M$,  denote the $i^{th}$-column of $M$ by \TCP{$M^{(i)}$}.
 First, $P\subseteq R^N$, is image of an idempotent matrix $\iota:R^N\lra R^N$. 
In particular, $P$ is generated by the columns $\iota^{(i)}$ of $\iota$. 
Denote $Q=\ker(\iota)$. We display two exact sequences and commutative diagrams of maps:
$$
\diagram
0\ar[r] & Q \ar[r] & R^N\ar[rd]_{\iota} \ar[r] & P\ar@{_(->}[d] \ar[r] & 0,\\
&&&R^N&\\
\enddiagram
~~
\diagram
0\ar[r] & Q^* \ar[r] & (R^N)^*\ar[dr]_{\iota^*} \ar[r] & P^*\ar[r]\ar@{^(->}[d] & 0\\
&&&(R^N)^*&\\
\enddiagram
$$
$$
{\rm Write}\quad R^N=\oplus_{i=1}^NRe_i, \quad (R^N)^*=\oplus_{i=1}^NRe_i^*.
$$
The latter diagram shows that $P^*$ is generated by the columns of $\iota^*$.
The quadratic structure on $P$ is given by 
an isomorphism $\varphi: P \lra P^*$. 
For $\sigma\in End(P, P^*)$, we extend $\sigma:R^N \lra (R^N)^*$, by defining $ \sigma_{|Q}=0$. 
So, we have 
$$
End(P, P^*) \subseteq End(R^N, (R^N)^*)=\BM_N(R). 
$$
So, $\sigma$ is given by a matrix $\Sigma\in \BM_N(R)$. Then, $\Sigma$ has the following properties:
(1)  $\Sigma^t=\Sigma$, 
(2) $\Sigma q=0$ for all $q\in Q$. 
(3) $\Sigma$ is injective on $P$,
(4) $image(\Sigma)=P^*$.

Let $\varepsilon_i:=(\iota^*)^{(i)}$  the columns of $\iota^*$. Then, $\{\varepsilon_i\}$ generates  $P^*$. Now,
\begin{equation}\label{incluNov10}
\exists ~p_i\in P~\ni~\varepsilon_i=\Sigma(p_i). 
\qquad {\rm Write}\quad p_i=\sum_j \lambda_{ij}(T)\iota^{(i)}
\end{equation}

Now, let $S$ be the set of all coefficients of entires in $\iota$,  $\Sigma$, $q$ and of $\lambda_{ij}$. Let $\BF$ be the prime field of $k$. 
Let $A_0=\BF[S]\subseteq A$ and $R_0=A_0[T_1, \ldots, T_n]$. Let $\iota_0:R_0^N\lra R_0^N$ be given by the matrix of $\iota$.
Let $P_0=image(\iota_0)$, $Q_0=\ker(\iota_0)$. Also, consider $\Sigma_0=\Sigma:R_0^n \lra (R_0^N)^*$. We have the commutative diagrams
$$
\diagram
&&&R_0^N&\\
0\ar[r] & Q_0\ar@{^(->}[d] \ar[r] & R_0^N\ar[ru]^{\iota_0} \ar[r]\ar@{^(->}[d]  & P_0\ar@{^(->}[d] \ar@{_(->}[u]\ar[r] & 0\\
0\ar[r] & Q \ar[r] & R^N\ar[rd]_{\iota} \ar[r] & P\ar@{^(->}[d] \ar[r] & 0\\
&&&R^N&\\
\enddiagram
\qquad 
\diagram
&&&(R_0^N)^*&\\
0\ar[r] & Q_0^* \ar[r] \ar@{^(->}[d]& (R_0^N)^*\ar@{^(->}[d]\ar[ru]^{\iota_0^*} \ar[r] & P_0^*\ar[r]\ar@{^(->}[d] \ar@{_(->}[u]& 0\\
0\ar[r] & Q^* \ar[r] & (R^N)^*\ar[dr]_{\iota^*} \ar[r] & P^*\ar[r]\ar@{^(->}[d] & 0\\
&&&(R^N)^*&\\
\enddiagram
$$
Now, Then,
\bE
\item From equation (\ref{incluNov10}), $\varepsilon_i\in P_0^*$ and $p_i\in P_0$.
\item $\Sigma_0^t=\Sigma^t=\Sigma=\Sigma_0$.
\item $\Sigma_0 q=0$ for all $q\in Q_0$. This is because $Q_0\subseteq Q$.
\item $\Sigma$ is injective on $P_0$. This is because $P_0\subseteq P$.
\item $image(\Sigma_0)=P_0^*$.\\
To see this, note from the diagram $Q_0=Q\cap R_0^N$. Likewise $P\cap R_0^N=P_0$ and
  $P^*\cap (R_0^N)^*=P_0^*$. So,
$$
\Sigma_0(R_0^N) \subseteq \Sigma(R^N)\cap (R_0^N)^*\subseteq P^*\cap (R_0^N)^*=P_0^*
$$
Note, $P_0^*$ is generated by the columns of $\iota^*$ (i. e. $\varepsilon_i$), as an $R_0$-module.
By equation (\ref{incluNov10}), since $p_i\in P_0$, we have $\Sigma_0$ maps on to $P_0^*$.

\eE
So, $\Sigma_0$ defines an isomorphism $\varphi_0: P_0\iso P_0^*$, such that $\varphi_0^*=\varphi_0$, because $\Sigma=\Sigma_0$
is symmetric. 
Also, $\varphi_0\otimes R=\varphi$, because they are restriction of the maps defined by $\Sigma_0=\Sigma$.
Therefore, \TCP{$(P_0, \varphi_0)\otimes R=(P, \varphi)$}. 

Now, we will enlarge $A_0$ to accommodate the local triviality condition.
There are $s_1, \ldots, s_m\in A$ and isometries $\psi_i: (A_{s_i}, q) \iso (\overline{P}_{s_i}, \varphi_{s_i})$ such that 
$s_1+\cdots+ s_m=1$.  Let $e_1, \ldots, e_r$ be the standard basis of $A^r$. 
Fix $l$ and work with $\psi_l: (A_{s_l}^r, q) \iso (\overline{P}_{s_l}, \varphi_{s_l})$. The generators of 
$\overline{P}$ is given by the columns $\overline{\iota^{(j)}}$. So, 
$$
\psi_l(e_i)=\sum_j\frac{a_{ijl}}{s_1^{\mu}}\overline{\iota^{(j)}}. \quad {\rm Also,}\quad 
\overline{\iota^{(i)}}= \sum_j\frac{b_{ijl}}{s_1^{\mu}}\psi_l(e_j)
$$
Write $A_1=A_0[s_i, a_{ijl}, b_{ijl}]$ and replace $A_0$ by $A_1$. Then, $\psi_l(e_i)\in \overline{P_0}_{s_l}$.
This defines a map $\eta_l: (A_0^l)_{s_i}\lra  \overline{P_0}_{s_l}$, by setting \TCP{$\eta_l(e_i) =\psi_l(e_i)$}. 
Then,   $\eta_l\otimes 1=\psi_l$.
By the latter equation on $\overline{\iota^{(i)}}$, it follows $\eta_i$ is surjective, and hence an isomorphism. 
To prove that $\eta_i$ is isometry, we need to check the commutativity of the diagrams
$$
\diagram
((A_0^r)_{s_l}, q)\ar[r]^{\eta_l}\ar[d] & (\overline{P}_0, \varphi_0)_{s_l}\ar[d]\\
((A_0^r)_{s_l}, q)^* & (\overline{P}_0, \varphi_0)_{s_l}^*\ar[l]^{\eta_l^*}\
\enddiagram
$$
The commutativity follows from Lemma \ref{HomPQ}.
This establishes that $(P_0, \varphi_0)$ is locally trivial.


By the Theorem of Popescu \cite{P}, we have the following commutative  the diagram of rings and homomorphisms 
$$
\diagram 
&& A_2\ar[d]\\
\BF \ar[r] & A_0 \ar[ru]^{\beta}\ar@{^(->}[r] & A\\
\enddiagram
\quad {\rm where}~A_2~{\rm is~smooth~over}~\BF.
$$
Write $R_2=A_2[T_1, \ldots, T_n]$. Then, by 
 (\ref{orthLindel}),  $(P_0, \varphi_0) \otimes R_2$ is extended from $A_2$. Since $(P, \varphi)= \left((P_0, \varphi_0) \otimes R_2\right)\otimes R$,
 we have $(P, \varphi)$ 
 is extended from $A$. This completes the proof.  
$\eop$


\section{Homotopy}\label{secHomotopy}

In this section, we use the Quadratic analogue (\ref{orthPopescu}) of Lindel's theorem to prove some key homotopy theorems.
First, we recall the following standard lemma. 
\bL\label{zLocTrivial}
Suppose $(A, \m)$ is a commutative noetherian local ring and ${\bf u}\in Q_{2n}'(A)$. Consider the orthogonal complement $K=A{\bf u}^{\perp}\subseteq (A^{2n+1}, q_{2n+1})$. 
Then, $K\cong (A, q_{2n})$,
where $q_{2n}=\sum_{i=1}^nX_iY_i$.
\eL
\pf Write ${\bf u}_0=(1; 0, \ldots, 0; 0, \ldots, 0)\in Q_{2n}'(A)$. 
By Lemma \ref{transSEO}, there is a matrix $\sigma\in {\SE}O(A, q_{2n+1})$,
such that ${\bf u}_0\sigma= {\bf u}$. So, we have a diagram of exact sequences
$$
\diagram 
0\ar[r] & K\ar[r] \ar@{-->}[d]_{\sigma_0} & A^{2n+1} \ar[r]^{\langle {\bf u}, -\rangle}\ar[d]_{\sigma}^{\wr} & A\ar[r]\ar@{=}[d] & 0\\
0\ar[r] & (A, q_{2n})\ar[r] & A^{2n+1} \ar[r]_{\langle {\bf u}_0, -\rangle} & A\ar[r] & 0\\
\enddiagram
$$
So, $\sigma$ induces an isometry $\sigma_0:K \iso (A, q_{2n})$. 
\pic $\eop$

\bT\label{PopThm3p4M2} 
Let $A$ be a regular ring over a  field $k$, with $1/2\in k$. 
Suppose $H(T)\in Q_{2n}'(A[T])$. Then, there is an orthogonal matrix $\sigma(T)\in O(A[T], q_{2n+1})$, such that 
$$
H(T)=H(0)\sigma(T) \qquad {\rm and}\qquad \sigma(0)=1.
$$
%
\eT
\pf 
Suppose $H(T)\in Q_{2n}'(A[T])$ is a homotopy. For $R= A[T], A$, use the following generic notations, to denote the quadratic modules
$$
\left\{
\begin{array}{lll}
q:=q_{2n+1}:R^{2n+1} \to R &{\rm sending} & (u_1, \ldots, u_n, v_1, \ldots, v_n, s)\mapsto \sum_{i=1}^nu_iv_i+s^2\\
q_0:R \to R &{\rm sending} & s\mapsto s^2\\
\end{array}
\right.
$$ 
As usual, define 
$B_q(e, e')=\frac{q(e+e')-q(e)-q(e')}{2}$. With respect to the standard basis, the matrix of $B_q$ is given by 
$$
B_q:=\frac{1}{2}\left(
\begin{array}{ccc} 
0 & I_n & 0\\ I_n & 0 & 0 \\ 0 & 0 & 2\\
\end{array}
\right)
$$
So, the map 
$
R^{2n+1} \to (R^{2n+1})^* ~{\rm sends} ~{\bf v} \mapsto {\bf v}B_q$. 
These bilinear forms give the following  exact sequences:
$$
\diagram
0\ar[r] & K \ar[r] & A[T]^{2n+1} \ar[rr]^{\langle H(T), -\rangle} & & A[T] \ar[r] & 0\\ 
0\ar[r] & K_0 \ar[r] & A^{2n+1} \ar[rr]_{\langle H(0), -\rangle} & & A \ar[r] & 0\\ 
\enddiagram
$$
Therefore, 
$K=\left(A[T]H(T)\right)^{\perp}$,  $K_0=\left(AH(0)\right)^{\perp}$
are orthogonal complements. 
Write $\overline{K}:=K\otimes \frac{A[T]}{(T)}$. By Lemma \ref{zLocTrivial}, $\overline{K}$ is locally isometric to $(A, q_{2n})$. By Theorem \ref{orthPopescu}, 
there is an isometry $\tau:K \iso \overline{K}\otimes A[T]$. 
Further, it follows 
$\overline{K}=\left(RH(0)\right)^{\perp}\cong K_0$. 
Therefore, there is an isometry $\sigma_0:  \overline{K} \iso K_0$, which extends to an isometry
$\sigma_0\otimes 1: \overline{K}\otimes A[T] \iso K_0\otimes A[T]$.
Then,  $\sigma_1:=(\sigma_0\otimes 1)\tau: K \iso  K_0\otimes A[T]$ is an isometry.
Finally, note 
$$
(A[T]H(T), q_{|A[T]H(T)}) \cong (A[T], q_0) \cong (A[T]H(0), q_{|A[T]H(0)}).
$$
Now, consider the  diagram
\begin{equation}\label{linearMaps}
\diagram
0\ar[r] & K\ar[r]\ar[d]_{\sigma_1} & A[T]^{2n+1}\ar@{-->}[d]^{\sigma(T)} \ar[rr]^{\langle H(T), -\rangle} & & A[T] \ar@{=}[d]\ar[r] & 0\\ 
0\ar[r] & K_0\otimes A[T] \ar[r] & A[T]^{2n+1} \ar[rr]_{\langle H(0), -\rangle} & & A[T] \ar[r] & 0\\ 
\enddiagram
\end{equation}
of quadratic spaces. In this diagram, the horizontal lines are split exact sequences of quadratic spaces. Hence, there is an 
isometry $\sigma(T)\in O\left(A[T], q\right)$, such that the diagram (\ref{linearMaps}) commutes.
Therefore $H(0)\sigma(T)=H(T)$.
 By construction, ({\it alternately, by replacing} $\sigma(T)$ by $\sigma(0)^{-1}\sigma(T)$), we have $\sigma(0)=1$.
\pic 
 $\eop$ 

The following Corollary would be of some use for our future discussions.
\bC\label{directHomo}
Let $A$ be a
regular ring over a   field  $k$, with $1/2\in k$ and $n\geq 2$ be an integer. 
Let ${\bf u}, {\bf v}\in Q_{2n}'(A)$ such that $[{\bf u}]= [{\bf v}]\in \pi_0(Q_{2n}')(A)$.
Then, there is a homotopy $H(T)\in Q_{2n}'(A[T])$ such that $H(0)={\bf u}$ and $H(1)={\bf v}$.
Equivalently, for 
${\bf u}, {\bf v}\in Q_{2n}(A)$ if $\zeta_0({\bf u})= \zeta_0({\bf v})\in \pi_0(Q_{2n})(A)$, then
there is a homotopy $H(T)\in Q_{2n}(A[T])$ such that $H(0)={\bf u}$ and $H(1)={\bf v}$.
\eC
\pf Suppose ${\bf u}, {\bf v}\in Q_{2n}'(A)$ such that $[{\bf u}]= [{\bf v}]\in \pi_0(Q_{2n}')(A)$.
Then, there is a sequence of homotopies $H_1(T), \ldots , H_m(T) \in Q_{2n}'(A[T], q_{2n+1})$
such that ${\bf u}_0:=H_1(0)={\bf u}$, ${\bf u}_m:=H_m(1)={\bf v}$ and $\forall~i=1, \ldots, m-1$, we have ${\bf u}_i:=H_i(1)=H_{i+1}(0)$.
By Theorem \ref{PopThm3p4M2}, for $i=1, \ldots, m$ there are orthogonal matrices $\sigma_i(T)\in O(A[T], q_{2n+1})$ such that
$\sigma_i(0)=1$ and 
$H_i(T)=H_i(0)\sigma_i(T) ={\bf u}_{i-1}\sigma_i(T)$. Therefore, ${\bf u}_i=H_i(1) = {\bf u}_{i-1}\sigma_i(1)$.
Write $H(T)={\bf u}_0\sigma_1(T)\cdots \sigma_m(T)$. Then, $H(T)\in Q_{2n}'(A[T])$ and $H(0)= {\bf u}_0$ and $H(1)= {\bf u}_m$.
This establishes first part of the statement on $\pi_0(Q_{2n}')(A)$. The latter assertion on $\pi_0(Q_{2n})(A)$ follows from the former, 
by the bijective correspondences 
$Q_{2n}'(A) \llra Q_{2n}(A)$ and $Q_{2n}'(A[T]) \llra Q_{2n}(A[T])$. 
\tcp $\eop$
\begin{remark}{\rm
Another way to state (\ref{directHomo}) would be to say that the homotopy relation on $Q_{2n}(A)$ is actually an equivalence relation.
}
\end{remark}

In a slightly more formal language, the above is summarized as follows.
\bT\label{summHomT}
Suppose $A$ is a
regular ring over a   field  $k$, with $1/2\in k$ and $n\geq 2$ is an integer. 
For, $\sigma(T)\in O(A[T], q_{2n+1})$ and ${\bf u}\in Q_{2n}'(A)$, define the  {\rm (right)} action 
${\bf u}\sigma(T):={\bf u}\sigma(1)\in Q_{2n}'(A)$. 
Denote $O(A, q_{2n+1}, T)=\{\sigma(T)\in O(A[T], q_{2n+1}): \sigma(0)=1 \}$.
Then, the map
$$
\frac{Q_{2n}'(A)}{O(A, q_{2n+1}, T)} 
\lra \pi_0(Q_{2n})(A) \quad {\rm is~a~bijection}.
$$
\eT
\pf Similar to the proof of (\ref{directHomo}). $\eop$
\begin{remark}{\rm 
A version of Theorem \ref{summHomT} is embedded in the proof of \cite[Theorem 1.0.6]{F}, where it was assumed that $A$ essentially 
smooth over an infinite field $k$, with $1/2\in k$. In \cite{F}, methods  of Nisnevich topology was used \cite{AHW}, 
while our methods are fairly basic.
}
\end{remark}

\section{The Involution}\label{secInvolution}
 
To begin with, we introduce the following key definition, in this article. 
\bD\label{invoDefi}
Suppose $A$ is a commutative ring. For\\
 ${\bf v}=(s; f_1, \ldots, f_n; g_1, \ldots, g_n)\in Q_{2n}(A)$ define
 $$
 \Gamma(s; f_1, \ldots, f_n; g_1, \ldots, g_n)= (1-s; f_1, \ldots, f_n; g_1, \ldots, g_n)
 $$
 This association, ${\bf v}\mapsto \Gamma({\bf v})$, establishes a bijective correspondence
 $$
 \Gamma: Q_{2n}(A)\iso Q_{2n}(A), \quad{\rm such ~that}\quad \Gamma^2=1_{Q_{2n}(A)}.
 $$
 That means, $\Gamma$ is an involution on $Q_{2n}(A)$. {\rm (This notation $\Gamma$ will be among the standard notations throughout this article.)}
\eD
We record the following  obvious lemma.
\bL\label{obvious}
Suppose $A$ is a commutative noetherian ring and $n\geq 2$ is an integer. Let 
 $\Gamma: Q_{2n}(A)\iso Q_{2n}(A)$ be the involution, as in (\ref{invoDefi})  and
 ${\bf v}=(s; f_1, \ldots, f_n; g_1, \ldots, g_n)\in Q_{2n}(A)$.
Then,
\bE
\item $\BI({\bf v})\cap  \BJ({\bf v})=(f_1, f_2, \ldots, f_n)$.
\item 
 $\BJ({\bf v})=\BI(\Gamma({\bf v}))$ and $\eta(\Gamma({\bf v}))=\left(\BJ({\bf v}), \omega_{\Gamma({\bf v})}\right)$.
\item For $H(T)\in Q_{2n}(A[T])$, we have
 $\Gamma(H(T))_{T=t}= \Gamma(H(t))$.
\item\label{obNov1216} Therefore,
$
\forall~{\bf v}, {\bf w}\in Q_{2n}(S)\qquad  \zeta_0({\bf v}) =\zeta_0({\bf w}) \Llra  \zeta_0(\Gamma({\bf v})) =\zeta_0(\Gamma({\bf w})).
$
\eE
\eL
We also record the following obvious observation.
\bL\label{Bibekananda}
Suppose $A$ is a commutative noetherian ring with $\dim A=d$. Suppose $n\geq 2$ is an integer.
Let ${\bf 0}\in Q_{2n}(A)$ be the base point and denote ${\bf 1}:=(\Gamma({\bf 0}))=(1; 0, \ldots 0; 0, \ldots, 0)$.
Then, $\zeta_0({\bf 0})= \zeta_0({\bf 1})$.
\eL
\pf For simplicity, we give the proof for $n=2$, which we exhibit in the following steps:
 \bE
 \item 
By Homotopy $H_1(T)=(0; T, 0; 0, 0)\in Q_4(A[T])$, we have \\
$\zeta_0(0; 0, 0; 0, 0)=\zeta_0(0; 1, 0; 0,  0)$.
 \item 
By Homotopy $H_2(T)=(T; 1, 0; T(1-T), 0)\in Q_4(A[T])$, we have\\
 $\zeta_0(0; 1, 0; 0,  0)=\zeta_0(1; 1,  0; 0, 0)$.

 \item 
By Homotopy $H_3(T)=(1; 1-T, 0; 0, 0)\in Q_4(A[T])$, we have\\
 $\zeta_0(1; 1, 0; 0,  0)=\zeta_0(1; 0, 0; 0,  0)$.
\eE
\pic $\eop$


\vspace{2mm}
In deed, $\Gamma$ factors through an involution on $\pi_0(Q_{2n})(A)$, as follows.
 
\bC\label{theInvolution} 
Suppose $A$ is a commutative ring and $n\geq 2$ is an integer. Then, the involution $\Gamma: Q_{2n}(A)\iso Q_{2n}(A)$ induces a bijective map
$\widetilde{\Gamma}: \pi_0(Q_{2n})(A) \iso \pi_0(Q_{2n})(A)$, which is also an involution {\rm (meaning, $\widetilde{\Gamma}^2=1_{\pi_0(Q_{2n})(A)}$)}.
{\rm (The notation $\widetilde{\Gamma}$ will also be among our standard notations throughout this article.)}
\eC 
\pf It follows from Lemma \ref{obvious}, that if $H(T)\in Q_{2n}(A(T))$ is a homotopy from ${\bf v}$ to ${\bf w}$, then $\Gamma(H(T))$ is a 
homotopy $\Gamma({\bf v})$ to $\Gamma({\bf w})$. The following commutative diagram of pushout squares,
 $$
 \diagram 
  Q_{2n}(A[T])\ar[dr]_{\Gamma} \ar[rr]^{T=0}\ar[dd]_{T=1} && Q_{2n}(A)\ar[dd]^{\zeta_0}\ar[dr]^{\Gamma}&\\
 &   Q_{2n}(A[T]) \ar[rr]^{T=0\qquad\qquad}\ar[dd]^{T=1} && Q_{2n}(A)\ar[dd]^{\zeta_0}\\
Q_{2n}(A) \ar[rr]^{\zeta_0\qquad\qquad}\ar[dr]_{\Gamma} && \pi_0\left(Q_{2n}\right)(A)\ar@{-->}[dr]^{\widetilde{\Gamma}}&\\ 
& Q_{2n}(A) \ar[rr]_{\zeta_0} && \pi_0\left(Q_{2n}\right)(A)\\ 
\enddiagram 
$$
in $\underline{Sets}$, establishes that $\widetilde{\Gamma}$ is well defined. Clearly, $\widetilde{\Gamma}^2=1$.
\pic $\eop$

\bC\label{altF_i4involution}
Suppose $A$ is a commutative ring and $n\geq 2$ is an integer. Suppose $(I, \omega)\in {\mathcal LO}(A, n)$.
Assume $f_1, \ldots, f_n\in I$ induce $\omega: A^n \sur \frac{I}{I^2}$ and $(f_1, \ldots, f_n)=I\cap J$,
with $I+J=A$.
Then, 
$$
\widetilde{\Gamma}(\zeta(I, \omega))= \zeta(J, \omega_J)\in \pi_0(Q_{2n})(A).
$$ 
\eC
 \pf Pick an element $s\in I$ such that $1-s\in J$. Then, $s(1-s)\in IJ=(f_1, \ldots, f_n)$. Also, $I=(f_1, \ldots, f_n, s)$
 and $J=(f_1, \ldots, f_n, 1-s)$. Write $s(1-s)=f_1g_1+\cdots+f_ng_n$. 
 With ${\bf v}=(s, f_1, \ldots, f_n, g_1, \ldots, g_n)\in Q_{2n}(A)$, we have 
 $$
 \widetilde{\Gamma}(\zeta(I, \omega))= \widetilde{\Gamma}(\zeta_0({\bf v})) = \zeta_0(\Gamma({\bf v})) =\zeta(J, \omega_J).
 $$
 \pic $\eop$
 
\begin{remark}\label{clarifyGamma}{\rm
We clarify a point that might get lost in  the formalisms of (\ref{theInvolution}) and (\ref{altF_i4involution}). Suppose
$(I, \omega_I)\in {\mathcal LO}(A, n)$ and ${\bf u}, {\bf v}\in Q_{2n}(A)$, with $\eta({\bf u})=\eta({\bf v})=(I, \omega_I)$.
We have $\zeta_0({\bf u})=\zeta_0({\bf v})=\zeta(I, \omega_I)$ (see \cite{F}). So, there is a sequence of 
homotopies $H_1(T), \ldots, H_m(T)\in Q_{2n}(A[T])$, such that, with $H_i(1)={\bf u}_i$,
we have  $H_1(0)={\bf u}=:{\bf u}_0$, $H_{i}(1)={\bf u}_i=H_{i+1}(0)$ and $H_m(1)={\bf u}_m= {\bf v}$. Now, 
$\Gamma(H_i(T))$ give a sequence of homotopies, starting from $\Gamma({\bf u})$ to $\Gamma({\bf v})$. This means, 
in the statement of (\ref{altF_i4involution}), $\zeta(J, \omega_J)$ depends on $(I, \omega_I)$ only, and is independent of
choice of $(J, \omega_J)$.
}
\end{remark}


The following is another version of the Moving Lemma \ref{movingLemm}. 
\bL[Moving Representation]\label{movingRep}
Suppose $A$ is a commutative noetherian ring with $\dim A=d$. 
Suppose $n\geq 2$ is an integer such that $2n\geq d+1$.
Let $K\subseteq A$ be an ideal with $height(K)\geq n$ and $x\in \pi_0(Q_{2n})(A)$. Then, there 
is a local $n$-orientation $(J, \omega_J)\in {\mathcal LO}(A, n)$ such that $x=\zeta(J, \omega_J)$, $height(J)\geq n$ and $J+K=A$.
\eL
\pf It is proved by two successive applications of Moving Lemma \ref{movingLemm}.
First, let $x=\zeta(I, \omega_I)$. By  (\ref{movingLemm}), there is 
${\bf u}=(s; f_1, \ldots, f_n; g_1, \ldots, g_n)$
such that $\eta({\bf u})=(I, \omega_I)$. Denote   $(I_0, \omega_{I_0}):=\eta(\Gamma({\bf u}))=(\BJ({\bf u}, \omega_{\Gamma({\bf u})})$. Then, 
$\tilde{\Gamma}(x)=\zeta(I_0, \omega_{I_0})$.

Now, we apply Moving Lemma \ref{movingLemm}, to $(I_0, \omega_{I_0})$ and $K$. 
There is ${\bf v}= (S; F_1, \ldots, F_n; G_1, \ldots, G_n)\in Q_{2n}(A)$, such that $\eta({\bf v})=(I_0, \omega_{I_0})$, and 
with $J=\BJ({\bf v})$, we have $height(J)\geq n$ and $J+K=A$. Now, $x=\tilde{\Gamma} (\tilde{\Gamma} (x))
= \tilde{\Gamma} (\zeta(I_0, \omega_{I_0}))=\zeta(J, \omega_J)$, where $\omega_J:=\omega_{\Gamma({\bf v})}$.
\pic $\eop$


We record the following useful lemma.
\bL\label{subsPri4vinQ2n}
Suppose $A$ is a commutative noetherian ring with $\dim A=d$. 
Suppose $n\geq 2$ is an integer.
Suppose ${\bf v}=(s; f_1, \ldots,f_n;  g_1, g_1, \ldots, g_n)\in Q_{2n}(A)$. 
Then, $\zeta_0({\bf v})=\zeta_0({\bf 0}) \Llra \zeta_0(\Gamma({\bf v}))=\zeta_0({\bf 0})$.

\eL
\pf By Lemma \ref{obvious} (\ref{obNov1216}), and Lemma \ref{Bibekananda}.
$$
\zeta_0({\bf v})=\zeta_0({\bf 0}) \Llra \zeta_0(\Gamma{\bf v})=\zeta_0({\bf 1})=\zeta_0({\bf 0}).
$$
 $\eop$

%


%
As a consequence of the above, we formulate 
  the following homotopy version of Subtraction Principle.
\bT\label{subPrincipeGamma} 
Suppose $A$ is a commutative noetherian ring with $\dim A=d$ and $n\geq 2$ is an integer.
Let $I, J$ be two ideals in $A$ such that
 $I+J=A$ and $IJ=(f_1, \ldots, f_n)$. 
 Let $(I, \omega_I), (J, \omega_J)\in {\mathcal LO}(A, n)$ be
  induced by $f_1, \ldots, f_n$.
  Then, $\zeta(I, \omega_I)) =\zeta_{0}({\bf 0}) \Llra \zeta(J, \omega_J)=\zeta_{0}({\bf 0})$.
\eT
\pf 
By local checking,  it follows that
$$
I=(f_1, \ldots, f_n)+I^2, \qquad{\rm and}\qquad J=(f_1, \ldots, f_n)+J^2.
$$
Since $I+J=A$, there is an element $s\in I$ such that $1-s\in J$. So, $s(1-s)\in IJ$ and hence 
$$
\sum_{i=1}^nf_ig_i+s(s-1)=0 \quad {\rm for~some} ~g_1, \ldots, g_n\in A.
$$
Therefore,
$$
 {\bf v}:=(s; f_1, \ldots, f_n; g_1, \ldots, g_n)\in Q_{2n}(A).
$$
We have $\zeta_0({\bf v}) =\zeta(I, \omega_I)$ and $\zeta_0(\Gamma({\bf v})) =\zeta(J, \omega_J)$.
Now, the proof follows from (\ref{subsPri4vinQ2n}).
$\eop$

\section{The Group Structure on $\pi_0(Q_{2n})(A)$} \label{secGroupStru}

In this section, we  establish a group structure on the set $\pi_0(Q_{2n})(A)$, when $2n\geq \dim A+2$ and $A$ is a regular ring 
over a field $k$, with $1/2\in k$.
We start with the following basic 
ingredient of the group structure. 
\bD\label{psudoPlus}
Let $A$ be a commutative noetherian ring and $n\geq 2$ be an integer.
{\rm (Refer to notations $\eta, \zeta_0, \zeta$ in  diagram \ref{zetaDiag}.)}
 Let
$(I, \omega_I), (J, \omega_J)\in {\mathcal LO}(A, n)$ be such that $I+J=A$. 
Let $\omega:=\omega_I\star \omega_J:  A^n \sur \frac{IJ}{(IJ)^2}$ be the unique surjective map induced by $\omega_I, \omega_J$.
We define a pseudo-sum 
$$
(I, \omega_I) \hat{+} (J, \omega_J) :=  \zeta(IJ, \omega)\in \pi_0(Q_{2n})(A).
$$
Also, for ${\bf u}, {\bf v}\in Q_{2n}(A)$ with $\BI({\bf u})+\BI({\bf v})=A$, define pseudo-sum
$$
{\bf u}\hat{+} {\bf v} := \eta({\bf u}) \hat{+} \eta({\bf v}) \in \pi_0(Q_{2n})(A).
$$
\eD 
The rest of this section is devoted to establish that this pseudo sum respects homotopy and extends to $\pi_0(Q_{2n})(A)$,
when 
$A$ is a regular ring 
over a  field $k$, with $1/2\in k$. 
The following is an obvious corollary.
\bC
Let $A$ be a commutative noetherian ring with $\dim A=d$ and $n\geq 2$ be an integer, \TCP{with $2n\geq d+3$.}
Suppose $(I, \omega_I), (J, \omega_J)\in \CL\CO(A, n)$ and $I+J=A$.
Let ${\bf u}=(s; f_1, \ldots, f_n; g_1, \ldots, g_n)\in Q_{2n}(A)$ and 
${\bf v}=(S; F_1, \ldots, F_n; G_1, \ldots, G_n)\in Q_{2n}(A)$ be such that $\eta({\bf u}) = (I, \omega_I)$
and $\eta({\bf v}) = (J, \omega_J)$. Write $\eta(\Gamma({\bf u}))= (I_1, \omega_{I_1}))$ and 
$\eta(\Gamma({\bf v}))= (J_1, \omega_{J_1}))$. We have $I+I_1=J+J_1=A$, and further assume 
$J+I_1=I+J_1=A$. {\rm ({\it Such choices of ${\bf u}, {\bf v}$ would be available, because  $2n\geq d+1$.})}
%
Then,
$$
\tilde{\Gamma}\left((I, \omega_I) \hat{+} (J, \omega_J) \right)=
(I_1, \omega_{I_1}) \hat{+} (J_1, \omega_{J_1}).
$$
\eC
\pf From addition principle (\cite[Theorem 5.6]{BK})
$$
I\cap I_1 \cap J \cap J_1=(U_1, \ldots, U_n)\quad \ni \quad  U_i-f_i\in (II_1)^2,~~ U_i-F_i\in (JJ_1)^2.
$$
Let  $\omega:A^n\sur \frac{IJ}{(IJ)^2}$, $\omega':A^n\sur \frac{I_1J_1}{(I_1J_1)^2}$ be the surjective maps
induced by $U_1, \ldots, U_n$. If follows $\omega= \omega_I\star \omega_J$ and $\omega'= \omega_{I_1}\star \omega_{J_1}$.
So, by corollary \ref{altF_i4involution},
 $$
 \tilde{\Gamma}\left((I, \omega_I) \hat{+} (J, \omega_J) \right)=
 \tilde{\Gamma}\left(\zeta(IJ, \omega)\right)=\zeta(I_1J_1, \omega')
 =
 (I_1, \omega_{I_1}) \hat{+} (J_1, \omega_{J_1}).
 $$
\pic $\eop$

\vspace{2mm}
\noindent{\bf Remark.} Note that the involution operations $\Gamma: Q_{2n}(A) \lra Q_{2n}(A)$ does not factor through a
map $\CL\CO(A, n)\lra \CL\CO(A, n)$.

\vspace{2mm} 
Now we define a pseudo-difference in the spirit of (\ref{psudoPlus}).
\bD\label{psudoDiff}
Suppose $A$ is a commutative ring and $n\geq 2$ be an integer.
Suppose $(K, \omega_K), (I, \omega_I) \in {\mathcal LO}(A, n)$. %
Assume that there is ${\bf u}=(s; f_1, \ldots, f_n; g_1, \ldots, g_n)\in Q_{2n}(A)$ such that 
$\eta({\bf u})=(I, \omega_I)$ and with $J=\BJ({\bf u})$, $J+K=A$. Therefore,
$$
(f_1, \ldots, f_n)=I \cap J \quad \ni\quad I+J=K+J=A. 
$$
Let $\omega_J:A^n\sur \frac{J}{J^2}$ be the surjective map induced by $f_1, \ldots, f_n$. 
Then, define the pseudo-difference
$$
(K, \omega_K) \hat{-} (I, \omega_I):=(K, \omega_K) \hat{+} (J, \omega_J) \in \pi_0(Q_{2n})(A).
$$
{\rm 
We remark:
(1)
A priori, the pseudo-difference depends on the choice of $J$.  There is no conditions on $height(I)$, $height(J)$, $height(K)$; nor did we assume $I+K=A$.
(2) By Moving Lemma \ref{movingLemm}, such choices  ${\bf u}\in Q_{2n}(A)$ would be available if $2n\geq \dim A+1$ and $height(K)\geq n$.

}
 \eD
 Subsequently, under additional hypotheses, we would first prove that the definition (\ref{psudoDiff})
 of pseudo difference does not depend of the choice of $(J, \omega_J)$. Then, we  
 prove that the pseudo difference is homotopy invariant with respect to either coordinate. A key to such proofs would be the following 
 lemma that  combines ({\it i.e. "adds"}) homotopies. ({\it For  basic element theory and the definition of generalized dimension functions,
 we refer to} \cite{M1}.)
\bL\label{HLK1hoK2}
Suppose $A$ is a commutative noetherian ring with $\dim A=d$ and $n\geq 2$ is an integer with $2n\geq d+2$. 
Consider a homotopy
 $$
 H(T)=(Z, \varphi_1, \ldots, \varphi_n, \gamma_1, \ldots, \gamma_n)\in Q_{2n}(A[T]).
 $$
  Write $\eta(H(0))=
(K_0, \omega_{K_0})$ and $\eta(H(1))=(K_1, \omega_{K_1})$.
Further suppose $(J, \omega_J)\in {\mathcal LO}(A, n)$ such that $K_0+J=K_1+J=A$ and $height(J)\geq n$. Then, there is a homotopy 
$\CH(T)\in Q_{2n}(A[T])$ such that $\eta(\CH(0))=(K_0J, \omega_{K_0J})$ and $\eta(\CH(1))=(K_1J, \omega_{K_1J})$,
where, for $i=0, 1$ $\omega_{K_iJ}:=\omega_{K_i}\star \omega_J:A^n \sur \frac{K_iJ}{(K_iJ)^2}$.
\eL
\pf  Write $Y=1-Z$. Then, $\BJ(H(T))= (\varphi_1, \ldots, \varphi_n, Y)$.
Write 
$$
\SP=\{\wp\in \spec{A[T]}: YT(1-T)\notin \wp,  J\subseteq \wp \}.
$$
 There is a generalized dimension function $\delta: \SP \lra \BN$ such that 
 $\forall~\wp\in \SP$,
 $\delta(\wp)\leq \dim\left(\frac{A[T]}{JA[T]}\right) \leq d+1-height(J) \leq d+1-n\leq n-1$.
 Further, $(\varphi_1, \ldots, \varphi_n, Y^2T(1-T))$ is a basic element in $A[T]^{n+1}$, on $\SP$. Therefore, there are polynomials 
$\lambda_1,\ldots, \lambda_n\in A[T]$ such that, with\\
 $\varphi_i'=\varphi_i+\lambda_iY^2T(1-T)$ for $i=1, \ldots, n$, 
we have 
$$
(\varphi_1', \varphi_2', \ldots, \varphi_n')
~~ {\rm is~basic~on}~\SP.
$$
We have \TCP{$\forall~i=1, \ldots, n~\varphi_i(0)= \varphi_i'(0),~\varphi_i(1)= \varphi_i'(1)$}.
 We compute 
$$
Z(1-Z)=Y(1-Y) =\sum_{i=1}^n \varphi_i\gamma_i= \sum_{i=1}^n \varphi_i'\gamma_i - Y^2T(1-T)\sum_{i=1}^n \lambda_i\gamma_i
$$
So,
$$
Y= \sum_{i=1}^n \varphi_i'\gamma_i  - Y^2T(1-T)\mu +Y^2\qquad {\rm where}\quad \mu = \sum_{i=1}^n \lambda_i\gamma_i.
$$
 Write $M=\frac{\BJ(H(T))}{(\varphi_1', \varphi_2', \ldots, \varphi_n')}$. Use "overline" to indicate images in 
$M$. We intend to repeat the proof of Nakayama's Lemma and we have 
$$
\left(\begin{array}{c}
\overline{\varphi_1}\\ \overline{\varphi_2}\\ \cdots \\ \overline{\varphi_n}\\   \overline{Y}\\
\end{array} \right)
=
\left(\begin{array}{ccccc}
0 & 0 & \cdots & 0 & -\lambda_1YT(1-T)\\
0 & 0 & \cdots & 0 & -\lambda_2YT(1-T)\\
\cdots &\cdots &\cdots &\cdots &\cdots\\
0 & 0 & \cdots & 0 & -\lambda_nYT(1-T)\\
0 & 0 & 0 & 0 & Y-  YT(1-T)\mu \\
\end{array} \right)
\left(\begin{array}{c}
\overline{\varphi_1}\\ \overline{\varphi_2}\\ \cdots \\ \overline{\varphi_n}\\ \overline{Y}\\
\end{array} \right)
\Lra 
$$
$$
\left(\begin{array}{ccccc}
1 & 0 & \cdots & 0 & \lambda_1YT(1-T)\\
0 & 1 & \cdots & 0 & \lambda_2YT(1-T)\\
\cdots &\cdots &\cdots &\cdots &\cdots\\
0 & 0 & \cdots & 1 & \lambda_nYT(1-T)\\
0 & 0 & 0 & 0 & 1-Y+ YT(1-T)\mu \\
\end{array} \right)
\left(\begin{array}{c}
\overline{\varphi_1}\\ \overline{\varphi_2}\\ \cdots \\ \overline{\varphi_n}\\ \overline{Y}\\
\end{array} \right)
=
\left(\begin{array}{c}
0 \\ 0 \\ \cdots \\ 0\\ 0 \\
\end{array} \right)
$$
Multiplying by the adjoint matrix and computing the determinant, with $Y'= Y-YT(1-T)\mu$, we have
$$
(1-Y')\BJ(H(T))\subseteq (\varphi_1', \varphi_2', \ldots, \varphi_n').
$$
We have \TCP{$Y'(0)=Y(0)=1-Z(0)$, $Y'(1)=Y(1)=1-Z(1)$.}
Further, 
$$
Y'(1-Y')=\sum_{i=1}^n\varphi'\gamma_i' \qquad{\rm for~some~polynomials}\quad \gamma_1', \ldots, \gamma_n'\in A[T]
$$
$$
H'(T)=(Y', \varphi_1', \ldots, \varphi_n'; \gamma'_1, \ldots, \gamma'_n) \in Q_{2n}(A[T]).
$$
We have
$$
\BJ(H(Y))= (\varphi_1, \ldots, \varphi_n, Y)= (\varphi_1', \ldots, \varphi_n', Y)=(\varphi_1', \ldots, \varphi_n', Y')=\BI(H'(T)).
$$
Claim
$$
\BJ(H'(T))+JA[T] =A[T].\qquad i.e. \qquad (\varphi_1', \ldots, \varphi_n', 1-Y')+JA[T] =A[T].
$$
To see this, let 
$$
\BJ(H'(T))+JA[T]\subseteq \wp \in \spec{A[T]}
$$
\bE
\item If $Y\in \wp$ then $(\varphi_1', \ldots, \varphi_n', Y)=(\varphi_1', \ldots, \varphi_n', Y')\subseteq \wp$. So, $Y'\in \wp$, which is 
impossible, since $1-Y'\in \wp$. \TCP{So, $\wp \in D(Y)$.}
\item 
Since $\wp \in D(Y)$  and since $(\varphi'_1, \ldots, \varphi'_n)$ is unimodular of $\SP$, we must have $T(1-T)\in \wp$.
\item Now, $T\in \wp$ implies, 
$$
\BJ(H'(0))+J= (\varphi_1'(0), \ldots, \varphi_n'(0), 1-Y'(0))+J
= (\varphi_1(0), \ldots, \varphi_n(0), 1-Y(0))+J
$$
$$
=(\varphi_1(0), \ldots, \varphi_n(0), Z(0))+J
=K_0+J=A\subseteq \wp,
$$
which is impossible. 
\item Likewise, $1-T\in \wp$ implies, 
$$
\BJ(H'(1))+J= (\varphi_1'(1), \ldots, \varphi_n'(1), 1-Y'(1))+J
= (\varphi_1(1), \ldots, \varphi_n(1), 1-Y(1))+J
$$
$$
=(\varphi_1(1), \ldots, \varphi_n(1), Z(1))+J
=K_1+J=A\subseteq \wp
$$
This is also impossible.
\eE
This establishes the claim.
So, 
$$
(\varphi_1', \ldots, \varphi_n')=\BI(H'(T))\cap \BJ(H'(T)) \quad {\rm where}\quad \BJ(H'(T))= (\varphi_1', \ldots, \varphi_n', 1-Y').
$$
Let 
$$
\Omega_{\BJ(H'(T))}:A[T]^n \sur \frac{\BJ(H'(T))}{\BJ(H'(T))^2} \quad {\rm be~induced~by}\quad \varphi_1', \ldots, \varphi'_n.
$$
Extend  $\omega_J:A^n \sur \frac{J}{J^2}$ to a surjective map  $\omega_{JA[T]}:A[T]^n \sur \frac{JA[T]}{J^2A[T]}$.
Let 
$$
\Omega:=\Omega_{\BJ(H'(T))}\star\omega_{JA[T]}  : A[T]^n \sur \frac{J\BJ(H'(T))}{J^2\BJ(H'(T))^2} \quad {\rm be~induced~by}~\Omega_{\BJ(H'(T))}, \omega_{JA[T]}.
$$
Now, $\Omega$ has a lift $\CH(T)\in Q_{2n}(A[T])$. Specializing at $T=0$ and $T=1$, we have
$$
\eta(\CH(0))= (K_0J, \omega_{K_0J}), \qquad \eta(\CH(1))= (K_1J, \omega_{K_1J}).
$$
The proof is complete. $\eop$


Now we proceed to prove, in several propositions, that the pseudo-difference (\ref{psudoDiff}) is well defined and homotopy invariant. 
\bP\label{wellpsDiff}
Suppose $A$ is a regular ring over a field $k$, with $1/2\in k$, with $\dim A=d$. 
Let $n\geq 2$ be an integer, such that $2n\geq d+2$. 
As in (\ref{psudoDiff}), let  $(K, \omega_K), (I, \omega_I)\in {\mathcal LO}(A, n)$ be given.

Let 
 ${\bf u}=(s, f_1, \ldots, f_n; g_1, \ldots, g_n)\in Q_{2n}(A)\in Q_{2n}(A)$ be such that, with 
$\eta({\bf u})=(I, \omega_I)$ and $\eta(\Gamma({\bf u}))=(J, \omega_J)$, $K+J=A$.
Likewise, let 
${\bf v}=(S, F_1, \ldots, F_n; G_1, \ldots, G_n)\in Q_{2n}(A)$ be such that
$\eta({\bf v})=(I, \omega_I)$, and with $\eta(\Gamma({\bf v}))=(L, \omega_L)$, $K+L=A$.
Assume  \TCP{$height(K)\geq n$}. 
Then, 
$$
(K, \omega_K) \hat{+} (J, \omega_J)= (K, \omega_K) \hat{+} (L, \omega_L) \in \pi_0(Q_{2n}(A).
$$
{\rm (Recall form (\ref{psudoDiff}), there is no restriction on $height(I), height(J)$ and, nor did we assume $I+K=A$.)}
\eP 

\pf 
By hypotheses,
$K+J=K+L=A$.
By Corollary \ref{altF_i4involution},   
$$
\widetilde{\Gamma}(\zeta(I, \omega_I))=[\zeta(J, \omega_J)]= [\zeta(L, \omega_L)]. ~~ {\rm Hence}~
\zeta_0(\Gamma({\bf u}))=\zeta_0(\Gamma({\bf v})). 
$$
By Corollary \ref{directHomo}, there is a homotopy
$$
H(T)=(Z(T); \varphi_1(T), \ldots, \varphi_n(T); \gamma_1(T), \ldots, \gamma_n(T)) \in Q_{2n}(A[T])
$$
such that $H(0)=\Gamma({\bf u})$, $H(1)=\Gamma({\bf v})$.
By the Homotopy Lemma \ref{HLK1hoK2} %
there is a homotopy $\CH(T)\in Q_{2n}(A[T])$ such that $\eta(\CH(0))=(KJ, \omega_K\star\omega_J)$ and $\eta(\CH(1))=(KL, \omega_K\star\omega_L)$.
This establishes,
$$
(K, \omega_K) \hat{+} (J, \omega_J)= (K, \omega_K) \hat{+} (L, \omega_L) \in \pi_0(Q_{2n}(A).
$$
\tcp $\eop$

\vspace{3mm}
\bC\label{wellPseDiff}
Suppose $A$ is a regular ring over a field $k$, with $1/2\in k$ and $\dim A=d$. 
Let $n\geq 2$ be an integer, such that $2n\geq d+2$.   Let $(K, \omega_K), (I, \omega_I) \in \CL\CO(A, n)$.
Assume  \TCP{$height(K)\geq n$}.
Then, the pseudo-difference {\rm (\ref{psudoDiff})}
$$
(K, \omega_K)\hat{-} (I, \omega_I) \in \pi_0(Q_{2n}(A)) \quad {\rm is~ well ~defined}.
$$
\eC
\pf This is immediate from Proposition \ref{wellpsDiff}. $\eop$


\vspace{3mm}

Now we prove that the pseudo-difference is homotopy invariant, with respect to the the $(I, \omega_I)$-coordiante. 
\bP\label{pseDinvHomo}
Suppose $A$ is a regular ring over a field $k$, with $1/2\in k$ and $\dim A=d$. 
Let $n\geq 2$ is an integer, with $2n\geq d+2$. 
Let\\
 $(K, \omega_K), (I_0, \omega_{I_0}), (I_1, \omega_{I_1}) \in {\mathcal LO}(A, n)$ and  \TCP{$height(K)\geq n$}.
Suppose\\
 $\zeta(I_0, \omega_{I_0})= \zeta(I_1, \omega_{I_1})\in \pi_0(Q_{2n}(A))$. Then, 
$$
(K, \omega_K) \hat{-} (I_0, \omega_{I_0})= (K, \omega_K) \hat{-} (I_1, \omega_{I_1}) \pi_0(Q_{2n}(A)).
$$
In other words, pseudo differences {\rm (defined in \ref{psudoDiff})} are homotopy invariant, with respect to the latter coordinate.
\eP
\pf By Moving Lemma \ref{movingLemm}, we can find ${\bf u}_0, {\bf u}_1\in Q_{2n}(A)$ such that $\eta({\bf u}_0)= (I_0, \omega_{I_0})$,
 $\eta({\bf u}_1)= (I_1, \omega_{I_1})$, and with $J_0=\BJ({\bf u}_0)$, $J_1=\BJ({\bf u}_1)$,  $K+J_0=K+J_1=A$. 
For $i=0, 1$, let $\omega_{J_i}=\omega_{\Gamma({\bf u}_i)}: A^n \sur \frac{J_i}{J_i^2}$. 
 By hypothesis,
 $\zeta(I_0, \omega_{I_0})= \zeta(I_1, \omega_{I_1})$ and hence
 $$
\zeta_0(\Gamma({\bf u}_0))= \zeta(J_0, \omega_{J_0})= \widetilde{\Gamma}(\zeta(I_0, \omega_{I_0})) =  \widetilde{\Gamma}(\zeta(I_1, \omega_{I_1})) =\zeta(J_1, \omega_{J_1})=\zeta_0(\Gamma({\bf u}_1)).
 $$
 By  Corollary \ref{directHomo}, there is a homotopy $H(T)\in Q_{2n}(A[T])$ such that $H(0)= \Gamma({\bf u}_0)$ and $H(1)=\Gamma({\bf u}_1)$. 
 Since $K+J_0=K+J_1=A$, by Lemma \ref{HLK1hoK2}, there is a homotopy $\CH(T)\in Q_{2n}(A[T])$ such that 
 $\eta(\CH(0))=(KJ_0, \omega_K\star\omega_{J_0})$ and $\eta(\CH(1))=(KJ_0, \omega_K\star\omega_{J_1})$.
Therefore,
$$
(K, \omega_K)\hat{+} (J_0, \omega_{J_0}) = (K, \omega_K)\hat{+} (J_1, \omega_{J_1})
$$
and hence
$$
(K, \omega_K)\hat{-} (I_0, \omega_{I_0}) = (K, \omega_K)\hat{-} (I_1, \omega_{I_1}).
$$
\pic$\eop$

\vspace{3mm}

We extend the definition of pseudo difference to $\pi_0(Q_{2n})(A)$ as follows.
\bC\label{psDiffKxpi0Q2nA}
Suppose $A$ is a regular ring over a field $k$, with $1/2\in k$, with $\dim A=d$. 
Let $n\geq 2$ is an integer, with $2n\geq d+2$. 
Let $(K, \omega_K) \in {\mathcal LO}(A, n)$ and  \TCP{$height(K)\geq n$}. Then, there is a well defined   set theoretic  map 
$$
\pi_0(Q_{2n})(A) \lra \pi_0(Q_{2n})(A) \qquad {\rm sending}\qquad \zeta(I, \omega_I)
\mapsto (K, \omega_K) \hat{-} (I, \omega_I).
$$
\eC
\pf Immediate from Proposition \ref{pseDinvHomo}. $\eop$

Now, we extend the pseudo-difference to $\pi_0(Q_{2n})(A) \times  \pi_0(Q_{2n})(A)$.

\bT\label{FINpsDiffThm}
Suppose $A$ is a regular ring over a field $k$, with $1/2\in k$ and  $\dim A=d$. 
Let $n\geq 2$ is an integer, with $2n\geq d+2$. Then, there is a set theoretic  map
$$
\Theta: \pi_0(Q_{2n})(A) \times  \pi_0(Q_{2n})(A) \lra  \pi_0(Q_{2n})(A)
$$
such that, for $(K, \omega_K) \in {\mathcal LO}(A, n)$, with  \TCP{$height(K)\geq n$}, and $(I, \omega_{I}) \in {\mathcal LO}(A, n)$,
\begin{equation}\label{DiaFINpsD}
\Theta\left(\zeta(K, \omega_K), \zeta(I, \omega_{I})\right)= (K, \omega_K)\hat{-}  (I, \omega_{I}).
\end{equation}
\eT
\pf Suppose $x\in \pi_0(Q_{2n})(A)$. 
By the Moving Lemma \ref{movingRep},
 we can write $x=\zeta(K, \omega_{K})$, with $height(K)\geq n$. Therefore, 
 if well defined,
 the Equation 
\ref{DiaFINpsD} applies to all $(x, y)\in \pi_0(Q_{2n})(A) \times  \pi_0(Q_{2n})(A)$.

Let $x, y\in \pi_0(Q_{2n})(A)$. We can write $x=\zeta(K, \omega_{K})$, with $height(K)\geq n$ and $y=\zeta(I, \omega_I)$.
Define,
$$
\Theta(x, y)= (K, \omega_K)\hat{-}  (I, \omega_{I}).
$$
We need to prove that, if $x=\zeta(K, \omega_{K})= \zeta(K', \omega_{K'})$, with $height(K)\geq n$ and $height(K')\geq n$,
then 
$$
(K, \omega_{K})\hat{-}  (I, \omega_{I})= (K', \omega_{K'})\hat{-}  (I, \omega_{I}).
$$
Again,  by Moving Lemma \ref{movingLemm}, there is ${\bf u}=(s; f_1, f_2, \ldots, f_n; g_1, \ldots, g_n) \in Q_{2n}(A)$ such that 
$\eta({\bf u})=(I, \omega_I)$, and with  $\eta(\Gamma({\bf u}))=(J, \omega_J)$,
we have $J+KK'=A$.
Since $x=\zeta(K, \omega_{K})=\zeta(K', \omega_{K.})$, %
by Corollary \ref{directHomo},  there is a homotopy $H(T)\in Q_{2n}(A[T])$ such that 
$\eta(H(0))=(K, \omega_K)$ and $\eta(H(1))=(K', \omega_{K'})$. By Lemma \ref{HLK1hoK2}, there is a homotopy 
$\CH(T)\in Q_{2n}(A[T])$ such that $\eta(\CH(0))=(KJ, \omega_K\star\omega_J)$ and $\eta(\CH(1))=(K'J, \omega_{K'}\star\omega_J)$. 
Therefore,
$$
(K, \omega_K)\hat{+}(J, \omega_J)=\zeta(KJ, \omega_K\star\omega_J)
= \zeta(K'J, \omega_{K'}\star\omega_J)= (K', \omega_{K'})\hat{+}(J, \omega_J).
$$
Therefore, by definition,
$$
(K, \omega_K)\hat{-}(I, \omega_I)= (K', \omega_{K'})\hat{-}(I, \omega_I).
$$
\tcp $\eop$



\vspace{3mm}

Finally, we are ready to define the group structure on $\pi_0(Q_{2n})(A)$.
\bD\label{defnBinary}
Suppose $A$ is a regular ring over a field $k$, with $1/2\in k$ and $\dim A=d$. 
Let $n\geq 2$ is an integer, with $2n\geq d+2$. 
Then, for $x, y\in \pi_0(Q_{2n})(A)$, the association
$$
(x, y) \mapsto \Theta\left(x, \widetilde{\Gamma}(y)\right)
$$ 
is a well defined binary operation on $\pi_0(Q_{2n})(A)$, where $\Theta$ is as in (\ref{FINpsDiffThm}). 
This operation is well defined because so are $\Theta$ and $\widetilde{\Gamma}$ (see Ccorollary \ref{theInvolution}).
We denote
$$
x+y:=\Theta\left(x, \widetilde{\Gamma}(y)\right)
$$
This binary operation will be referred to an addition.
\eD
With the help of the Moving Lemma \ref{movingLemm}, the addition operation on $\pi_0(Q_{2n})(A)$ can be described in a more direct manner,
as follows.
\bL\label{describeSum}
Suppose $A$ is a regular ring over a field $k$, with $1/2\in k$ and $\dim A=d$. 
Let $n\geq 2$ is an integer, with $2n\geq d+2$. 
Let $x, y\in \pi_0(Q_{2n})(A)$. By the Moving Lemma \ref{movingLemm}, $x=\zeta(K, \omega_K)$ and
$y=\zeta(I, \omega_I)$, for some  $(K, \omega_K), (I, \omega_I) \in {\mathcal LO}(A, n)$, such that $I+K=A$ and
 $height(K)\geq n$.
Then, 
$$
x+y= (K, \omega_K)\hat{+} (I, \omega_I) \qquad {\rm as~in~Definition~(\ref{psudoPlus})}.
$$
\eL
\pf 
Let ${\bf u}\in Q_{2n}(A)$ be such that $\eta({\bf u})=(I, \omega_I)$ and  write $\eta({\bf u})=(J, \omega_J)$.
Then, $\widetilde{\Gamma}(\zeta(I, \omega_I))= \zeta(J, \omega_J)$.
By Definition \ref{defnBinary},
$$
x+y=\Theta\left(x, \widetilde{\Gamma}(y)\right)= \Theta\left(x, \zeta(J, \omega_J)\right)= (K, \omega_K)\hat{-} (J, \omega_J)=  (K, \omega_K)\hat{+} (I, \omega_I).
$$
\tcp $\eop$


\vspace{2mm}
The following is a final statement on the group structure on $\pi_0(Q_{2n})(A)$.
\bT\label{abelianGroup}
Suppose $A$ is a regular ring over a field $k$, with $1/2\in k$ and $\dim A=d$. 
Let $n\geq 2$ is an integer, with $2n\geq d+2$. Then, the addition operation on $\pi_0(Q_{2n})(A)$, defined in (\ref{defnBinary}), endows a
structure of an abelian group on $\pi_0(Q_{2n})(A)$.
\eT
\pf First, ${\bf 0}=\zeta(A, \omega_A)=\zeta_0(0; 0, \ldots, 0; 0, \ldots, 0)=\zeta_0(1; 0, \ldots, 0; 0, \ldots, 0) $ acts as the additive identity of this addition. 
Given $x, y, z\in \pi_0(Q_{2n})(A)$, by applications of the Moving Lemma \ref{movingLemm}, we can write 
$$
x=\zeta(K, \omega_K),~y=\zeta(I, \omega_I)~z=\zeta(J, \omega_J)\quad \ni~K+I=K+J=I+J=A
$$
and $height(K)\geq n$, $height(I)\geq n$, $height(J)\geq n$. By Lemma \ref{describeSum}, we have the following.
$$
(x+y)+z=((K, \omega_K)\hat{+}(I, \omega_I))\hat{+} (J, \omega_J)=x+(y+z).
$$
So, the associativity holds. 
Further,
$$
x+y= (K, \omega_K)\hat{+}(I, \omega_I)= (I, \omega_I)\hat{+}(K, \omega_K)=y+x.
$$
So, the commutativity holds. Also,
$$
x+\widetilde{\Gamma}(x)= \Theta(x, x)= 0.
$$
So, $x$ has an additive inverse. \tcp $\eop$

\subsection{The Euler Class Groups}\label{secular}



In this subsection 
we compare the Euler class groups $E^n(A)$ and $\pi_0(Q_{2n})(A)$. 
Refer to Section \ref{eulerDefSec} for the definition of Euler class groups  and some notations used in this section.
First, we define a map $\rho: E^n(A) \lra \pi_0(Q_{2n})(A)$, as follows.

\bD\label{EntopizeroNov14}{\rm 
Suppose $A$ is a regular  ring over a field $k$, with $1/2\in k$ and  $\dim A=d$. Let $n\geq 2$ be  an integer, with $2n\geq d+2$.
Then, the  restriction of the map $\zeta: {\mathcal LO}(A, n) \lra \pi_0(Q_{2n})(A)$ to ${\mathcal LO}^c(A, n)$ ({\it see {\rm (Definition \ref{newEuler})} for this notation}) defines a 
 set theocratic map ${\mathcal LO}^c(A, n) \lra \pi_0(Q_{2n})(A)$. Since $\pi_0(Q_{2n})(A)$ has the structure of an abelian group (\ref{describeSum}), 
 this extends to a group homomorphism $\rho_0: \BZ\left({\mathcal LO}^c(A, n)\right)\lra \pi_0(Q_{2n})(A)$. Then,
  $\rho_0$ induces a surjective homomorphism $\rho: E^n(A)\sur\pi_0(Q_{2n})(A)$.
}
\eD 
\pf That $\rho_0$ factors through a homomorphism $\rho:E^n(A)\sur\pi_0(Q_{2n})(A)$ is obvious. 
For $x\in \pi_0(Q_{2n})(A)$, by Moving Lemma \ref{movingLemm}, $x=\zeta(I, \omega_I)$ for some 
$(I, \omega_I)\in {\mathcal LO}(A, n)$, with $height(I)\geq n$. It follows, by (\ref{abelianGroup}), $x=\rho(\overline{\varepsilon}(I, \omega_I))$.
\tcp $\eop$

\vspace{3mm}
\begin{remark}{\rm 
It follows from \cite[Theorem 4.2]{BS2}, if $\rho$  in (\ref{EntopizeroNov14}) 
is an isomorphism and $2n\geq d+3$, then  for $(I, \omega_I) \in {\mathcal LO}(A, n)$, with $height(I) \geq n$,
$\zeta(I, \omega_I) ={\bf 0}$ implies $\omega_I$ lifts to a surjective map $A^n\sur I$.

Conversely, the homomorphism $\rho$ in (\ref{EntopizeroNov14}) is clearly an isomorphism, if for $(I, \omega_I) \in {\mathcal LO}(A, n)$, and $height(I) \geq n$, 
$\zeta(I, \omega_I)=0$ implies $\omega_I$ lifts to a surjective map $A^n\sur I$. 
This would be case, in the following case.

}
\end{remark}

\bT\label{useBKNov14}
Suppose $k$ is an infinite perfect field, with $1/2\in k$ and $A$ is an essentially smooth ring over $k$, with $\dim A=d$. 
Let $n\geq 2$ be an integer, with $2n\geq d+3$. Then, the homomorphism $\rho$  in (\ref{EntopizeroNov14}) is an isomorphism. 
\eT
\pf We only need to prove that $\rho$ is injective. Let $\rho(x)=0$ for some $x\in E^n(A)$. We can write $x=\overline{\varepsilon}(I, \omega_I)$,
for some $(I, \omega_I) \in {\mathcal LO}(A, n)$, with $height(I)\geq n$. Let
 ${\bf u}\in Q_{2n}(A)$ be such that $\eta({\bf u})=(I, \omega_I)$. By (\ref{directHomo}), there is a homotopy $H(T)
 \in Q_{2n}(A[T])$ such that $H(0)={\bf 1}$ and $H(1)={\bf u}$. Write $H(T)=(Z; f_1, \ldots, f_n; g_1, \ldots, g_n)$. 
 By Moving Lemma, we can modify $f_1, \ldots, f_n$ by multiples of $T(1-T)Z^2$ and assume $height(\BI(H(T)))\geq n$.
 Write $\eta(H(T))=(\CI, \omega)$. Then, $\eta(H(0))=\eta({\bf 1})= (A, *)$, where $*$ is the trivial map. 
 Since $\CI+TA[T]=A[T]$, we have $\frac{\CI}{T\CI}\iso \CI(0)=A$. 
  Let $\lambda:A^n \sur A$ be the map given by $(1, 0, \ldots, 0)$. Then, $\omega$ and $\lambda$ combines to give give 
  surjective homomorphism $\varphi_0:A[T]^n \sur \frac{\CI}{T\CI^2}$. By \cite[Theorem 4.13]{BK}, $\varphi_0$ lifts to a 
  surjective map $\varphi:A[T]^n\sur \CI$. Specializing at $T=1$, we have $\varphi(1): A^n\sur \CI(0)=I$ is a lift of $\omega_I$.
Therefore $\overline{\varepsilon}(I, \omega_I)=0\in E^n(A)$. \tcp $\eop$

\appendix

\section{Preliminaries on Euler Class Groups} \label{eulerDefSec}
For a regular ring $A$ with $\dim A=d$, Nori outlined a definition of an alternate obstruction group $E^d(A)$, of codimension $d$-cycles,
which was expected to coincide with $\pi_0(Q_{2d})(A)$. Generators of the groups were $S=\{(m, \omega_m)\in {\mathcal LO}(A,d):
m\in \max(A) \}$. By some innovative analogies, for each co-dimension $n \geq 0$, Bhatwadekar and Sridharan \cite{BS2} defined an obstructions group,
 where $A$ was assumed to be any noetherian commutative ring. These groups have come to be known as 
 Euler Class Groups. 
In this section we modify the definition in \cite{BS2}, to avoid some superfluous elements.
({\it In the sequel, for a set $S$, the free abelian group generated by $S$ will be denoted by} $\BZ(S)$).
\bD\label{newEuler}
{\rm
Suppose $A$ is a noetherian commutative ring, with $\dim A=d$ and $n\geq 0$ is an integer. Let \\
${\mathcal LO}^c(A, n)=\{(I, \omega_I)\in {\mathcal LO}(A, n):  V(I)~{\rm is~connected~and}~height(I)\geq n \}$.

Suppose $(I, \omega_I)\in {\mathcal LO}(A, n)$ and $I=\cap_{i=1}^mI_i$ be a decomposition, where $V(I_i)\subseteq 
\spec{A}$ are connected. For $i=1, \ldots, m$, the orientation $(I, \omega_I)$ induce orientations $(I_i, \omega_{I_i})\in {\mathcal LO}(A, n)$.
Denote $\varepsilon(I, \omega_I)=\sum_{i=1}^m(I_i, \omega_{I_i})\in \BZ\left({\mathcal LO}^c(A, n)\right)$.
Recall,
an orientation $(I, \omega_I)\in {\mathcal LO}(A, n)$ is called global, if $\omega_I$ lifts to a surjective map $A^n\sur I$.
 Let $\SR^n(A)$ denote the subgroup of $\BZ\left({\mathcal LO}^c(A, n)\right)$, generated by the set \\
$\left\{\varepsilon(I, \omega_I):   (I, \omega_I)~{\rm is~global~and} ~height(I)\geq n\right\}$.
Define
$$
E^n(A)= \frac{\BZ\left({\mathcal LO}^c(A, n)\right)}{\SR^n(A)}
$$
Images of $\varepsilon(I, \omega_I)$ in $E^n(A)$ with be denoted by $\overline{\varepsilon}(I, \omega_I)$.

The definition of the Euler class group in \cite{BS2}, which we denote by $\CE^n(A)$,    differs a little from the above, which we recall.

Given an ideal $I$, two local orientations $(I, \omega), (I, \omega')\in {\mathcal LO}(A, n)$ were defined to be 
equivalent, if $\omega'=\omega\sigma$ for some elementary matrix $\sigma\in EL_n(A)$. The equivalence class of 
$(I, \omega)$ is denoted by $(I, [\omega])$. Let
${\mathcal LO}^c_0(A, n)=\{(I, [\omega_I]): (I, \omega_I)\in {\mathcal LO}(A, n):  V(I)~{\rm is~connected~and}~height(I)\geq n \}$.
Again, for $(I, \omega_I)\in {\mathcal LO}(A, n)$ and $I=\cap_{i=1}^mI_i$ be a decomposition, where $V(I_i)\subseteq 
\spec{A}$ are connected, as above, denote 
$\varepsilon_0(I, \omega_I)=\sum_{i=1}^m(I_i, [\omega_{I_i}])\in \BZ\left({\mathcal LO}^c_0(A, n)\right)$.
Let $\SR^n_0(A)$ denote the subgroup of $\BZ\left({\mathcal LO}_0^c(A, n)\right)$, generated by the set 
$\left\{\varepsilon_0(I, \omega_I):   (I, \omega_I)~{\rm is~global~and} ~height(I)\geq n\right\}$.
In \cite{BS2}, the Euler class group $\CE^n(A)$ was defined as
$$
\CE^n(A)= \frac{\BZ\left({\mathcal LO}_0^c(A, n)\right)}{\SR_0^n(A)}.
$$

}
\eD
However, the differences between $E^n(A)$ and $\CE^n(A)$ is superfluous, as follows.
\bL\label{superFlu}
Use notations as in (\ref{newEuler}). The map $\varphi: \BZ\left({\mathcal LO}^c(A, n)\right)\lra \BZ\left({\mathcal LO}_0^c(A, n)\right)$
sending $(I, \omega_I) \mapsto (I, [\omega_I])$ induces  
  a surjective homomorphism $\Phi_n:  E^n(A) \sur \CE^n(A)$. If $2n\geq d+3$, the $\Phi_n$ is an isomorphism.
\eL
\pf It is obvious that $\Phi_n$ is a well defined surjective homomorphism. To prove injectivity, let $\Phi_n(x)=0$.
Since, $2n\geq d+1$, by Moving Lemma \ref{movingLemm} $x=\overline{\varepsilon}(I, \omega_I)\in E^n(A)$, for some 
$(I, \omega_I)\in {\mathcal LO}(A, n)$, with $height(I)\geq n$. It follows, 
$image(\varepsilon_0(I, \omega_I))=0 \in \CE^n(A)$. By \cite[Theorem 4.2]{BS2}, it follows $\omega_I$ lifts to a surjective map 
$A^n\sur I$. Therefore, $(I, \omega_I)$ is global and $x=\overline{\varepsilon}(I, \omega_I)=0$. Hence, $\Phi_n$ is an isomorphism. 
$\eop$

\section{Elements of Elementary Orthogonal Subgroups}\label{EOsection}
This section is included due to nonavailability  of  suitable references.
In this section, for a commutative ring $A$,   the elements of $A^{n}$ would be considered as column
matrices. We consider the Orthogonal groups $O\left(A^{2n+1}, q_{2n+1}\right) \subseteq GL_{2n+1}(A)$ and $O\left(A^{2n}, q_{2n}\right)
 \subseteq GL_{2n}(A)$.
In analogy to Elementary Subgroups of $GL_n(A)$, Elementary Orthogonal subgroups ${\SE}O\left(A^{2n+1}, q_{2n+1}\right)$ and 
${\SE}O\left(A^{2n}, q_{2n}\right)$ of the respective Orthogonal  groups were defined,  which is classical (see \cite{CF}). 

For $i, j=0, 1, \ldots, 2n$, $i\neq j$ and $\lambda\in A$, let $E_{ij}(\lambda)\in GL_{2n+1}(A)$ denote the matrix, whose diagonal entries are $1$ 
and only other nonzero entry is $\lambda$ at the $(i, j)^{th}$ position. Recall, the Elementary subgroup $EL_{2n+1}(A)\subseteq  GL_{2n+1}(A)$
is generated by the set 
$$
\left\{E_{ij}(\lambda): i, j=0, 1, \ldots, 2n;~i\neq j; \lambda\in A \right\}
$$
For $i, =0, 1, \ldots, 2n$, $j=1, \ldots, 2n$, $i\neq j$   and $\lambda\in A$, one can write down an orthogonal matrices 
$\varepsilon_{ij}(\lambda) \in O\left(A^{2n+1}, q_{2n+1}\right)$, by modifying $E_{ij}(\lambda)$, in a natural way. The elementary 
orthogonal subgroup ${\SE}O\left(A^{2n+1}, q_{2n+1}\right)$ of
$O\left(A^{2n+1}, q_{2n+1}\right)$ would  generated by the set 
$$
\left\{\varepsilon_{ij}(\lambda): i=0, 1, \ldots, 2n; j= 1, \ldots, 2n;~i\neq j; \lambda\in A \right\},
$$
while some $\varepsilon_{ij}(\lambda)$ would be redundant. A more precise definitions of $\varepsilon_{ij}(\lambda)$ are as follows.
\bD\label{fiveGenDef}{\rm
Suppose $A$ is a commutative ring and $n\geq 1$ is an integer. 
For our purpose, elements of $A^{2n+1}$ would be dented by the columns 
$$
{\bf u}:=(z; x_1, \ldots, x_n; y_1, \ldots, y_n)^t.
$$
For  convenience, these $2n+1$ coordinates would be referred to as $0^{th}, 1^{st}, 2^{nd}, \ldots, n^{th}, \ldots, (2n)^{th}$-coordinates. 

\bE
\item \label{TypeOne} For $j=1, \ldots, n$, let $\varepsilon_{0, j}(\lambda)$ be the matrix that sends 
$$
{\bf u}\mapsto (z+\lambda x_i; x_1, \ldots, x_n; y_1, \ldots, y_{i-1},y_i-2\lambda z- \lambda^2x_i, y_{i+1},\ldots, y_n)^t
$$
For example, with $n=2$,  and $j=1$, we have 
$$
\varepsilon_{0, 1}(\lambda)=
\left( 
\begin{array}{ccccc}
1 & \lambda & 0 & 0 & 0\\
0 & 1 & 0 & 0 & 0\\
0 & 0 & 1 & 0 & 0\\
-2\lambda &  -\lambda^2 & 0 & 1 & 0\\
0 & 0 & 0 & 0 & 1\\
\end{array}
\right)
$$
\item \label{TypeTwo}  For $j=n+1, \ldots, 2n$, let $\varepsilon_{0, j}(\lambda)$ be the matrix that sends 
$$
{\bf u}\mapsto (z+\lambda y_i; x_1, \ldots, x_{i-1}, x_i-2\lambda z- \lambda^2y_i, x_{i+1},\ldots, x_n; y_1, \ldots, y_n)^t
$$
For example, with $n=2$,  and $j=n+1=3$, we have 
$$
\varepsilon_{0, 3}(\lambda)=
\left( 
\begin{array}{ccccc}
1 & 0 & 0 & \lambda & 0\\
-2\lambda & 1 &0 & -\lambda^2  & 0\\
0 & 0 & 1 & 0 & 0\\
0 &  0 & 0 & 1 & 0\\
0 & 0 & 0 & 0 & 1\\
\end{array}
\right)
$$
\item \label{TypeThree}  For $i, j=1, \ldots, n$, $i\neq  j$, let
$\varepsilon_{i, j}(\lambda)$ be the matrix that sends 
$$
{\bf u}\mapsto (z; x_1, \ldots, ,x_{i-1}, x_i+\lambda x_j , x_{i+1}, \ldots, x_n; y_1, \ldots, y_{j-1}, y_j-\lambda y_i, y_{j+1}, \ldots,  y_n )^t
$$
For example, with $n=2$,  and $i=1, j=2$, we have 
$$
\varepsilon_{1, 2}(\lambda)=
\left( 
\begin{array}{ccccc}
1 & 0 & 0 & 0 & 0\\
0 & 1 &\lambda & 0  & 0\\
0 & 0  & 1 & 0 & 0\\
0 &  0 & 0 & 1 & 0\\
0 & 0 & 0& -\lambda & 1\\
\end{array}
\right)
$$
\item\label{TypeFour}  For $i=1, \ldots, n$; $j=n+1, \ldots, 2n$; $i\neq j-n$  ({\it without loss of generality} $i < j-n$), let
$\varepsilon_{i, j}(\lambda)$ be the matrix that sends 
$$
{\bf u}\mapsto (z; x_1, \ldots,   x_i+\lambda y_j ,  \ldots, x_{j-n} -\lambda y_i, \ldots, x_n; y_1, \ldots,  y_n )^t
$$ 
For example, with $n=2$,  and $i=1, j=n+2=4$, we have 
$$
\varepsilon_{1, 4}(\lambda)=
\left( 
\begin{array}{ccccc}
1 & 0 & 0 & 0 & 0\\
0 & 1 &0 & 0  & \lambda\\
0 & 0  & 1 & -\lambda & 0\\
0 &  0 & 0 & 1 & 0\\
0 & 0 & 0& 0 & 1\\
\end{array}
\right)
$$
\item \label{TypeFive}  For $i=n+1, \ldots, 2n$; $j=1, \ldots, n$; $i-n\neq j$  ({\it without loss of generality} $i-n < j$), let
$\varepsilon_{i, j}(\lambda)$ be the matrix that sends 
$$
{\bf u}\mapsto (z; x_1, \ldots,  x_n; y_1, \ldots,   y_{i-n}+\lambda x_j ,  \ldots, y_j -\lambda x_{i-n}, \ldots, y_n )^t
$$ 
For example, with $n=2$,  and $i=n+1=3, j=2$, we have 
$$
\varepsilon_{3, 2}(\lambda)=
\left( 
\begin{array}{ccccc}
1 & 0 & 0 & 0 & 0\\
0 & 1 &0 & 0  & 0\\
0 & 0  & 1 & 0  & 0\\
0 &  0 & \lambda & 1 & 0\\
0 & -\lambda & 0& 0 & 1\\
\end{array}
\right)
$$
\eE
}
\eD
Now we are ready to give the sedition of the Elementary Orthogonal subgroup ${\SE}O\left(A, q_{2n+1}\right) \subseteq O\left(A, q_{2n+1}\right)$.

\bD\label{EODEF}
Suppose $A$ is a commutative ring. Then, the Elementary Orthogonal subgroup ${\SE}O\left(A, q_{2n+1}\right)$ of the Orthogonal 
group $O\left(A, q_{2n+1}\right)$ is defined to be the subgroup, generated by the
set of all $\varepsilon_{i, j}(\lambda)$ described and listed above {\rm (\ref{fiveGenDef})}.

{\rm Likewise, the Elementary Orthogonal subgroup ${\SE}O\left(A, q_{2n}\right)$ of the Orthogonal 
group $O\left(A, q_{2n}\right)$ is defined.}
\eD

We give a proof of the following well known lemma, which we use in the proof of 
Lemma \ref{zLocTrivial}.
\bL\label{transSEO}
Suppose $(A, \m)$ is a local commutative ring, with $1/2\in A$. Let
$$
{\bf u}:=(z; x_1, \ldots, x_n; y_1, \ldots, y_n)^t\in Q_{2n}'(A)\quad {\rm and}\quad 
{\bf u}_0:=(1; 0, \ldots, 0; 0, \ldots, 0)^t.
$$
Then, there is a matrix $\sigma\in {\SE}O\left(A, q_{2n+1}\right)$ such that $\sigma {\bf u}={\bf u}_0$.
\eL
\pf We prove this in several step. 

\noindent{\bf Step 1}. We prove that there is a matrix $\tau\in {\SE}O\left(A, q_{2n+1}\right)$, such that 
$\tau {\bf u}=(s; ~a_1, \ldots, a_n;~b_1, \ldots, b_n)^t$, with $a_1$ unit.

If $x_1$ is a unit, there is nothing to prove. So, assume $x_1\in \m$.
We have $z^2+\sum_{i=1}^nx_iy_i=1$. Since $x_1\in \m$, at least one of $z, x_2y_2, \ldots, x_ny_n$ would be a unit.
\bE
\item Suppose   $y_j$ is unit, for some $j=2, \ldots, n$. Without loss go generality, assume $y_2$ is a unit. 
Then, (see Type \ref{TypeFour}) $\tau:=\varepsilon_{1,n+2}(1)$ has the desired property.
\item Now, assume $y_j\in \m, ~\forall~j=2, \ldots, n$ and $y_1$ is a unit. Then, (see Type \ref{TypeThree})
$\varepsilon_{1, 2}(1)({\bf u})=(z, x_1+x_2, x_2, \ldots, x_n;~y_1, y_2-y_1, y_3, \ldots, y_n)^t$. Then, $y_2':=y_2-y_1$ is a unit.
If $x_1':=x_1+x_2$ is a unit, then $\tau:=\varepsilon_{1, 2}(1)$ would have the desired property. If $x_1'\in \m$ then, $\tau= \varepsilon_{1,n+2}(1)\varepsilon_{1, 2}(1)$
would have the desired property.
\item So, we assume $y_j\in \m, ~\forall~j=1, \ldots, n$  and hence $z$ is a unit. 
Then, $\varepsilon_{0, n+1}(1)({\bf u})=(z+y_1; x_1-2z-y_1; x_2, \ldots, x_n;~y_1, \ldots, y_n)$. Note, $x_1':=x_1-2z-y_1$ is unit. So,
$\tau=\varepsilon_{0, n+1}(1)$ has the desired property. 

\eE
This completes the proof of Step 1.

\noindent{\bf Step 2.} We prove that there is a matrix $\tau\in {\SE}O\left(A, q_{2n+1}\right)$, such that 
$\tau {\bf u}=(1; ~a_1, \ldots, a_n;~b_1, \ldots, b_n)^t$, with $a_1$ unit.

To see this note that, by Step 1, we can assume that $x_1$ is a unit. Now, with $\lambda=x_1^{-1}(1-z)$,
let $\tau:=\varepsilon_{0, 1}(\lambda)$. Then, $\tau({\bf u}=(1, x_1, \ldots, x_n;~y_1-2z\lambda-\lambda^2x_1, y_2, \ldots, y_n)^t$.
So, $\tau$ has the desired property of Step 2.

\noindent{\bf Step 3.} We prove that there is a matrix $\tau\in {\SE}O\left(A, q_{2n+1}\right)$, such that 
$\tau {\bf u}=(1; ~a_1,0, \ldots,  \ldots, 0;~b_1, \ldots, b_n)^t$, with $a_1$ unit.

By Step 3, we can assume ${\bf u}=(1, x_1, x_2, \ldots, x_n;~y_1, \ldots, y_n)^t$, with $x_1$ unit.
For $i=2, \ldots, n$, write (see Type \ref{TypeThree})  $\tau_i=\varepsilon_{1,i}\left(-x_1^{-1}x_i\right)$ and let $\tau=\tau_n\tau_{n-1}\cdots \tau_2$.
Then $\tau({\bf u})= (1, x_1, 0, \ldots, 0;~y_1', y_2, \ldots, y_n)^t$ for some $y_1'$. This established Step 3.

\noindent{\bf Step 4.} We prove that there is a matrix $\tau\in {\SE}O\left(A, q_{2n+1}\right)$, such that 
$\tau {\bf u}=(1; ~a_1,0, \ldots,  \ldots, 0;~0, \ldots, 0)^t$, with $a_1$ unit.

By Step 3, we can assume ${\bf u}=(1,~ x_1, 0, \ldots, 0;~y_1, y_2, \ldots, y_n)^t$, with $x_1$ unit. Since ${\bf u}\in Q_{2n}'(A)$,
we have $1+x_1y_1=1$ and hence $y_1=0$.
Therefore, ${\bf u}=(1,~ x_1, 0, \ldots, 0;~0, y_2, \ldots, y_n)^t$, with $x_1$ unit. 
For $i=2, \ldots, n$, write (see Type \ref{TypeFive})  $\tau_i=\varepsilon_{n+i, 1}\left(-x_1^{-1}y_i\right)$ and let $\tau=\tau_n\tau_{n-1}\cdots \tau_2$.
Since $x_2=x_3=\cdots=x_n=0$ and $y_1=0$, we have 
$\tau({\bf u})= (1,~ x_1, 0, \ldots, 0;~0, 0, \ldots, 0)^t$,with $x_1$ unit.
This establishes Step 4.

\noindent{\bf The Final Step.} By Step 4, we can assume ${\bf u}= (1; ~x_1,0, \ldots,  \ldots, 0;~0, \ldots, 0)^t$, with $x_1$ unit.
Now, let $\tau=\varepsilon_{0, n+1}\left(\frac{x_1}{2}\right)$. Then, $\tau({\bf u})=(1; 0, \ldots, 0;~0, \ldots, 0)$.

\pic $\eop$.

\end{document}